\documentclass[letterpaper, 10 pt, conference]{ieeeconf}  

\IEEEoverridecommandlockouts                              

\usepackage{amsmath}
\usepackage{amssymb}
\usepackage{tikz}
\usepackage{pgfplots}
\usepackage{graphicx}
\usepackage{setspace}
\usepackage{caption}
\usepackage{subcaption}
\usepackage{xcolor}

\newcommand{\until}[2]{\, \mathcal{U}_{[{#1},{#2}]} \, }

\newcommand{\rs}[1]{\rho^{#1}(\boldsymbol{x},t)}

\newcommand{\rsq}[1]{\rho^{#1}(\boldsymbol{x},t_q)}
\newcommand{\rsoo}[1]{\rho^{#1}(\boldsymbol{x},0)}
\newcommand{\rsnt}[1]{\rho^{#1}(\boldsymbol{x})}
\newcommand{\rsnts}[1]{\rho^{#1}(\boldsymbol{x}(t_*))}
\newcommand{\rsntt}[1]{\rho^{#1}\big(\boldsymbol{x}(t)\big)}
\newcommand{\rsnto}[1]{\rho^{#1}(\boldsymbol{x}_0)}
\newcommand{\rss}[2]{\rho^{#1}(\boldsymbol{x},{#2})}
\newcommand{\rso}[1]{\rho^{#1}_{opt}}

\newcommand{\re}[1]{{\color{black}#1}}

\newtheorem{definition}{Definition} 
\newtheorem{theorem}{Theorem} 
\newtheorem{assumption}{Assumption} 
\newtheorem{problem}{Problem} 
\newtheorem{remark}{Remark}
\newtheorem{lemma}{Lemma}

\newtheorem{example}{Example}

\title{\LARGE \bf
Prescribed Performance Control \\ for Signal Temporal Logic Specifications
}

\author{Lars Lindemann, Christos K. Verginis, and Dimos V. Dimarogonas
\thanks{This work was supported in part by the Swedish Research Council (VR), the European Research Council (ERC), the Swedish Foundation for Strategic Research (SSF),  the EU H2020 Co4Robots project, and the Knut and Alice Wallenberg Foundation (KAW).}
\thanks{The authors are with the Department of Automatic Control, School of Electrical Engineering, Royal Institute of Technology (KTH), 100 44 Stockholm, Sweden. {\tt\small llindem@kth.se (L. Lindemann), cverginis@kth.se (C.K. Verginis), dimos@kth.se (D.V. Dimarogonas)}}%
}

\begin{document}

\maketitle
\thispagestyle{empty}
\pagestyle{empty}

\begin{abstract}

Motivated by the recent interest in formal methods-based control for dynamic robots, we discuss the applicability of prescribed performance control to nonlinear systems subject to signal temporal logic specifications. Prescribed performance control imposes a desired transient behavior on the system trajectories that is leveraged to satisfy atomic signal temporal logic specifications. A hybrid control strategy is then used to satisfy a finite set of these atomic specifications. Simulations of a multi-agent system, using consensus dynamics, show that a wide range of specifications, i.e., formation, sequencing, and dispersion, can be robustly satisfied.

\end{abstract}

\section{Introduction}
\label{sec:introduction}

Temporal logics have lately gained much attention in robotic applications due to the possibility of formulating complex temporal specifications leading to formal methods-based control strategies \cite{belta2007symbolic,fainekos2009temporal}. These logics have for instance been used in multi-agent systems to perform realistic real-world tasks such as sequencing, coverage, surveillance, and formation control. In this multi-agent setup, linear temporal logic (LTL) \cite{guo2017task,kloetzer2010automatic} and metric interval temporal logic (MITL) \cite{nikou2016cooperative} have been used. These approaches abstract the physical environment, including robot dynamics, and the temporal logic formula into a finite-state automaton representing all possible robot motions. Search algorithms are then used to find a formula-satisfying discrete path that is subsequently accomplished by continuous control laws. However, these approaches may be subject to the state-space explosion problem \cite[Section 2.3]{baier}.

Robustness of temporal logic formulas was discussed in \cite{fainekos2009robustness} with the introduction of the robustness degree and the robust semantics, which are an under-approximation of the robustness degree. These notions give a measure of how robustly a formula is satisfied, i.e., a continuous scale indicating if a formula is marginally or greatly satisfied. Signal temporal logic (STL) \cite{maler2004monitoring} uses quantitative time properties and entails space robustness \cite{donze2}, a form of robust semantics. 

Prescribed performance control (PPC) \cite{bechlioulis2008robust,bechlioulis2014low} explicitly takes the transient and steady-state behavior of a tracking error into account. A user-defined performance function prescribes a desired temporal behavior that is then achieved by a continuous state feedback control law.

STL was introduced in the context of monitoring \cite{maler2004monitoring,donze2}, but not control. Control of systems subject to STL is a difficult task due to the nonlinear, nonconvex, noncausal, and nonsmooth semantics. Previous work on STL control synthesis has been done in \cite{lindemann,raman1,sadraddini} by using model predictive control (MPC), while \cite{liu2017communication}  explicitly extends the method derived in \cite{raman1} to multi-agent systems. In this paper, we consider a nonlinear system subject to a subset of STL. We propose to recast this constrained control problem into a PPC framework to satisfy atomic temporal formulas. Subsequently, the hybrid system framework in \cite{goebel2012hybrid,goebel2009hybrid} is used to satisfy a finite set of these atomic temporal formulas. To the best of the authors' knowledge, the approach presented in this paper is the first approach using a continuous state feedback control law for STL specifications. 

The remainder of this paper is organized as follows: Section \ref{sec:backgound} introduces notation and preliminaries. Section \ref{sec:problem_formulation} illustrates the underlying main idea and the problem definition. Section \ref{sec:strategy} presents a control law satisfying atomic temporal formulas, while Section \ref{sec:sequencing} considers a finite set of these atomic temporal formulas. Section \ref{sec:simulations} presents simulations of a centralized multi-agent system subject to different STL formulas, followed by a conclusion in Section \ref{sec:conclusion}.
This an extended version of the 56th IEEE Conference on Decision and Control (2017) version.
\section{Notation and Preliminaries}
\label{sec:backgound}

Scalars are denoted by lowercase, non-bold letters $x$ and column vectors are lowercase, bold letters $\boldsymbol{x}$. The vector $\boldsymbol{0}_n$ consists of $n$ zeros. True and false are denoted by $\top$ and $\bot$ with $\mathbb{B}:=\{\top,\bot\}$; $\mathbb{R}^n$ is the $n$-dimensional vector space over the real numbers $\mathbb{R}$. The natural, non-negative, and positive real numbers are $\mathbb{N}$, $\mathbb{R}_{\ge0}$, and $\mathbb{R}_{>0}$, respectively. 

\subsection{Signals and Systems}
Let $\boldsymbol{x}\in\mathbb{R}^n$, $\boldsymbol{u}\in\mathbb{R}^m$, and $\boldsymbol{w}\in \mathcal{W}$ be the state, input, and additive noise of a nonlinear system 
\begin{align}\label{system_noise}
\dot{\boldsymbol{x}}&=f(\boldsymbol{x})+g(\boldsymbol{x})\boldsymbol{u}+\boldsymbol{w},
\end{align}
where $\mathcal{W}\subset\mathbb{R}^n$ is a bounded set and the functions $f$ and $g$ satisfy Assumption \ref{assumption:1}.

 \begin{assumption}\label{assumption:1}
	The functions $f:\mathbb{R}^n\to\mathbb{R}^n$ and $g:\mathbb{R}^n\to\mathbb{R}^{n\times m}$ are  locally Lipschitz continuous, and $g(\boldsymbol{x})g^T(\boldsymbol{x})$ is positive definite for all $\boldsymbol{x}\in \mathbb{R}^n$.
\end{assumption} 
 
For the upcoming analysis, two basic results regarding the existence of solutions for initial-value problems (IVP) are needed. Assume $\boldsymbol{y}\in\Omega_{\boldsymbol{y}}\subseteq \mathbb{R}^{n+1}$ and consider the IVP
\begin{align}
\dot{\boldsymbol{y}} := H(\boldsymbol{y},t) \label{eq:xi_systems} \text{ with } \boldsymbol{y}(0):=\boldsymbol{y}_0\in \Omega_{\boldsymbol{y}},
\end{align}
where $H:\Omega_{\boldsymbol{y}}\times \mathbb{R}_{\ge 0} \to \mathbb{R}^{n+1}$ and $\Omega_{\boldsymbol{y}}$ is a non-empty and open set. A solution to this IVP is a signal $\boldsymbol{y}:\mathcal{J}\to\Omega_{\boldsymbol{y}}$ with $\mathcal{J}\subseteq \mathbb{R}_{\ge 0}$ obeying \eqref{eq:xi_systems}. In this paper, we will not explicitly distinguish between the state $\boldsymbol{y}$ and the solution $\boldsymbol{y}$ of \eqref{eq:xi_systems}.

\begin{lemma}\cite[Theorem 54]{sontag2013mathematical}\label{theorem:sontag1}
Consider the IVP in \eqref{eq:xi_systems}. Assume that $H:\Omega_{\boldsymbol{y}}\times \mathbb{R}_{\ge 0}\to \mathbb{R}^{n+1}$ is: 1) locally Lipschitz on $\boldsymbol{y}$ for each $t\in \mathbb{R}_{\ge 0}$, 2) piecewise continuous on $t$ for each fixed $\boldsymbol{y}\in\Omega_{\boldsymbol{y}}$.  Then, there exists a unique and maximal solution $\boldsymbol{y}:\mathcal{J}\to\Omega_{\boldsymbol{y}}$ with $\mathcal{J}:=[0,\tau_{max})\subseteq \mathbb{R}_{\ge 0}$ and $\tau_{max}\in \mathbb{R}_{>0}\cup\infty$.
\end{lemma}

\begin{lemma}\cite[Proposition C.3.6]{sontag2013mathematical}\label{theorem:sontag2}
Assume that the assumptions of Lemma \ref{theorem:sontag1} hold. For a maximal solution $\boldsymbol{y}$ on $\mathcal{J}=[0,\tau_{max})$ with $\tau_{max}<\infty$ and for any compact set $\Omega_{\boldsymbol{y}}^\prime\subset\Omega_{\boldsymbol{y}}$, there exists $t^\prime\in \mathcal{J}$ such that $\boldsymbol{y}(t^\prime)\notin\Omega_{\boldsymbol{y}}^\prime$.
\end{lemma}  

\subsection{Signal Temporal Logic (STL)}
Signal temporal logic is a predicate logic based on continuous-time signals. STL consists of predicates $\mu$ that are obtained after evaluation of a function $h:\mathbb{R}^n\to\mathbb{R}$ as 
$
 \mu:=
 \begin{cases} 
 \top \text{ if } h(\boldsymbol{x})\ge 0\\
 \bot \text{ if } h(\boldsymbol{x})< 0.
 \end{cases}
$ Note that $\boldsymbol{x}$ is seen here as a state and not a signal. For instance, consider the predicate $\mu:=(x\ge 1)$, which can be expressed by $h(x):=x-1$. Hence, $h$ determines the truth value of $\mu$ and maps from $\mathbb{R}^n$ to $\mathbb{R}$, while $\mu$ maps from $\mathbb{R}^n$ to $\mathbb{B}$. The STL syntax, given in Backus-Naur form, is 
\begin{align}
\phi \; ::= \; \top \; | \; \mu \; | \; \neg \phi \; | \; \phi_1 \wedge \phi_2 \; | \; \phi_1  \until{a}{b} \phi_2\;,
\end{align}
where $\mu$ is a predicate and $\phi_1$, $\phi_2$ are STL formulas. The temporal until-operator $\until{a}{b}$ is time bounded with time interval $[a,b]$ where $a,b\in \mathbb{R}_{\ge 0}\cup \infty$ such that $a\le b$.  The semantics of STL are introduced in Definition \ref{def:01} where the satisfaction relation $(\boldsymbol{x},t)\models \phi$ denotes that the signal $\boldsymbol{x}:\mathbb{R}_{\ge 0}\to\mathbb{R}^n$, possibly a solution of \eqref{system_noise} with $\boldsymbol{x}_0:=\boldsymbol{x}(0)$, satisfies  $\phi$ at time $t$.

\begin{definition} \cite[Definition 1]{maler2004monitoring} The STL semantics are recursively given by:
\begin{align*}
&(\boldsymbol{x},t) \models \mu 				 	&\Leftrightarrow	\;\;\; 	&h(\boldsymbol{x}(t))\ge 0\\
&(\boldsymbol{x},t) \models \neg\mu 			 	&\Leftrightarrow	\;\;\; 	&\neg((\boldsymbol{x},t) \models \mu)\\
&(\boldsymbol{x},t) \models \phi_1 \wedge \phi_2 	 	&\Leftrightarrow \;\;\;	 	&(\boldsymbol{x},t) \models \phi_1 \wedge (\boldsymbol{x},t) \models \phi_2\\
&(\boldsymbol{x},t) \models \phi_1 \until{a}{b} \phi_2		&\Leftrightarrow \;\;\;	 	&\exists t_1 \in[t+a,t+b]\text{ s.t. }(\boldsymbol{x},t_1)\models \phi_2 \\
&									&					&\wedge \forall t_2\in[t,t_1]\text{,}(\boldsymbol{x},t_2) \models \phi_1
\end{align*}
\label{def:01} 
\end{definition}
\vspace{-6mm}

The disjunction-, eventually-, and always-operator can be derived as $\phi_1\vee\phi_2=\neg(\neg\phi_1\wedge\neg\phi_2)$, $F_{[a,b]}\phi=\top  \until{a}{b} \phi$, and $G_{[a,b]}\phi = \neg F_{[a,b]}\neg\phi$. Additionally, robust semantics have been introduced in \cite{fainekos2009robustness} as a robustness measure. Space robustness \cite{donze2} $\rs{\phi}$ are robust semantics for STL given in Definition \ref{def:2}, for which it holds that $(\boldsymbol{x},t)\models \phi$ if $\rs{\phi}>0$.  Space robustness determines how robustly a signal $\boldsymbol{x}$ satisfies the formula $\phi$.
\begin{definition}\cite[Definition 3]{donze2} {The semantics of space robustness are recursively given by:}
\begin{align*}
\rs{\mu}& := h(\boldsymbol{x}(t))\\
\rs{\neg\phi} &:= 	-\rs{\phi}\\
\rs{\phi_1 \wedge \phi_2} &:= 	\min\big(\rs{\phi_1},\rs{\phi_2}\big)
\end{align*}
\begin{align*}
\rs{\phi_1 \vee \phi_2} &:= \max\big(\rs{\phi_1},\rs{\phi_2}\big)\\
\rs{\phi_1 \until{a}{b} \phi_2} &:= \underset{t_1\in[t+a,t+b]}{\max} \bigg( \min\Big(\rss{\phi_2}{t_1},\\ 
&\hspace{2.2cm} \underset{t_2\in[t,t_1)}{\min}\rss{\phi_1}{t_2} \Big) \bigg) \\
\rs{F_{[a,b]} \phi} &:= \underset{t_1\in[t+a,t+b]}{\max}\rss{\phi}{t_1}\\
\rs{G_{[a,b]} \phi} &:= \underset{t_1\in[t+a,t+b]}{\min}\rss{\phi}{t_1}.
\end{align*}
\label{def:2}
\end{definition}

 We abuse the notation as $\rho^\phi(\boldsymbol{x}(t)):=\rs{\phi}$ if $t$ is not explicitly contained in $\rs{\phi}$. For instance, $\rho^\mu(\boldsymbol{x}(t)):=\rs{\mu} := h(\boldsymbol{x}(t))$ since $h(\boldsymbol{x}(t))$ does not contain $t$ as an explicit parameter. However, $t$ is explicitly contained in $\rs{\phi}$ if temporal operators (eventually, always, or until) are used. In this paper, conjunctions are approximated by smooth functions.

\begin{assumption}\label{assumption:2}
The non-smooth conjunction $\rs{\phi_1 \wedge \phi_2}$ in Definition \ref{def:2} is approximated by a smooth function as $\rs{\phi_1 \wedge \phi_2}\approx -\ln\big(\exp(-\rs{\phi_1})+\exp(-\rs{\phi_2})\big)$.
\end{assumption}

\begin{remark}\label{remark:over_approximation}
The aforementioned approximation is an under-approximation of the robust semantics in Definition \ref{def:2}, i.e., $-\ln\big(\exp(-\rs{\phi_1})+\exp(-\rs{\phi_2})\big)\le \min\big(\rs{\phi_1},\rs{\phi_2}\big)$. This means that $(\boldsymbol{x},t) \models \phi_1 \wedge \phi_2$ if $-\ln\big(\exp(-\rs{\phi_1})+\exp(-\rs{\phi_2})\big)> 0$.
\end{remark}

\subsection{Prescribed Performance Control (PPC)}
Prescribed performance control (PPC) \cite{bechlioulis2008robust,bechlioulis2014low} constrains a generic error $\boldsymbol{e}:\mathbb{R}_{\ge 0}\to\mathbb{R}^n$ to a funnel. For instance, consider $\boldsymbol{e}(t)=\boldsymbol{x}(t)-\boldsymbol{x}_d(t)$ where $\boldsymbol{x}_d$ is a desired trajectory. In order to prescribe transient and steady-state behavior to this error, define the performance function $\gamma$ in Definition \ref{def:p}.
\begin{definition}\label{def:p}\cite{bechlioulis2014low} 
	A performance function $\gamma:\mathbb{R}_{\ge 0}\to\mathbb{R}_{> 0}$ is a continuously differentiable, bounded, positive, and non-increasing function. We define
	$\gamma(t):=(\gamma_0-\gamma_\infty)\exp(-lt)+\gamma_\infty$ where $\gamma_0, \gamma_\infty \in \mathbb{R}_{>0}$ with $\gamma_0\ge \gamma_\infty$ and $ l \in \mathbb{R}_{\ge 0}$.
\end{definition}

The task is to synthesize a feedback control law such that, given $-\gamma_i(0)<e_i(0)<M\gamma_i(0)$, the errors $e_i$ satisfy
\begin{align}
-\gamma_i(t)<e_i(t)<M\gamma_i(t) \;\;\;   \forall t\in \mathbb{R}_{\ge 0}, \forall i\in \{1,\hdots,n\}\label{eq:constrained_funnel}
\end{align}
with $0\le M\le 1$ and $\gamma_i$ as in Definition \ref{def:p}; $\gamma_i$ is a design parameter by which transient and steady-state behavior of $e_i$ can be prescribed. Similar to $M$ in the right inequality of \eqref{eq:constrained_funnel}, another constant could be added to the left inequality, which however will not be considered here. Note also that \eqref{eq:constrained_funnel} is a constrained control problem with $n$ constraints subject to the dynamics in \eqref{system_noise}. Next, define the normalized error $\xi_i:=\frac{e_i}{\gamma_i}$ and the transformation function $S$ as in Definition \ref{def:SS}. 
\begin{definition}\label{def:SS}
A transformation function $S:(-1,M)\to\mathbb{R}$ is a strictly increasing function, hence injective and admitting an inverse. In particular, we define $S(\xi):=\ln\left(-\frac{\xi+1}{\xi-M}\right)$.
\end{definition}

Dividing \eqref{eq:constrained_funnel} by $\gamma_i$ and applying the transformation function $S$ results in an unconstrained control problem  $-\infty<S\big(\xi_i(t)\big)<\infty$ with the transformed error $\epsilon_i:=S\big(\xi_i\big)$. If $\epsilon_i(t)$ is bounded for all $t$, then $e_i$ satisfies \eqref{eq:constrained_funnel}. This is a consequence of the fact that $S$ admits an inverse.

\section{Casting STL Control into a PPC framework}
\label{sec:problem_formulation}
We consider a subset of STL, which is expressive enough to formulate many real-world specifications. Considering the predicate $\mu$, the syntax is
\begin{subequations}\label{eq:subclass}
\begin{align}
\psi \; &::= \; \top \; | \; \mu \; | \; \neg \mu \; | \; \psi_1 \wedge \psi_2\label{eq:psi_class}\\
\phi \; &::= \;  G_{[a,b]}\psi \; | \; F_{[a,b]} \psi \label{eq:phi_class}\\
\theta^{s_1} &::= \bigwedge \limits_{i=1}^{N} \phi_i \text{ with } b_n\le a_{n+1}, \; \forall n\in\{1,\hdots,N-1\} \label{eq:theta1_class} \\
\theta^{s_2} &::=  F_{[c_1,d_1]} \Big( \psi_1 \wedge F_{[c_2,d_2]}\big(\psi_2 \wedge F_{[c_3,d_3]} (\hdots \wedge \phi_N) \big) \Big)\label{eq:theta2_class}\\
\theta &::= \theta^{s_1} \; |\; \theta^{s_2},\label{eq:theta_class}
\end{align}
\end{subequations}
where $\psi_1$ and $\psi_2$ are formulas of class $\psi$, whereas $\phi_i$ with $i\in\{1,\hdots,N\}$ are formulas of class $\phi$ with time intervals $[a_i,b_i]$. This STL subset is in positive normal form \cite{baier} and does not use disjunction- or until-operators. We refer to $\psi$ as non-temporal formulas. Due to the previous discussion, we write $\rsntt{\psi} :=\rs{\psi}$ and sometimes even omit $t$ resulting in $\rho^\psi(\boldsymbol{x})$. In contrast, $\phi$ and $\theta$ are referred to as temporal formulas due to the use of the always- and eventually-operators. We further refer to formulas \eqref{eq:phi_class} by the term atomic temporal formulas, while formulas in \eqref{eq:theta_class} are denoted as sequential formulas. Note that \eqref{eq:theta_class} either consists of \eqref{eq:theta1_class} or \eqref{eq:theta2_class}. 

\begin{assumption}\label{assumption:4}
Each formula of class $\psi$ that is contained in \eqref{eq:phi_class}, \eqref{eq:theta1_class}, and \eqref{eq:theta2_class} is: 1) s.t. $\rho^\psi(\boldsymbol{x})$ is concave and 2) well-posed in the sense that $(\boldsymbol{x},0)\models \psi$ implies $\|\boldsymbol{x}(0)\|<\infty$.
\end{assumption}
\begin{remark}
Part 2) of Assumption  \ref{assumption:4} is not restrictive since $\psi_{Ass.3}:=(\|\boldsymbol{x}\|<c)$, where $c$ is a sufficiently large positive constant, can be combined with the desired $\psi$ so that $\psi\wedge\psi_{Ass.3}$ is well-posed.
\end{remark}

The first objective in this paper is to synthesize a continuous feedback control law $\boldsymbol{u}(\boldsymbol{x},t)$ for atomic temporal formulas $\phi$ in \eqref{eq:phi_class} such that $\rsoo{\phi} > r$ where $r\in\mathbb{R}_{\ge 0}$ is a robustness measure and $\boldsymbol{x}:\mathbb{R}_{\ge 0}\to\mathbb{R}^n$ is the closed-loop solution of \eqref{system_noise} with initial condition $\boldsymbol{x}_0$. Additionally, we will upper bound $\rsoo{\phi}<\rho_{max}$ with $\rho_{max}\in\mathbb{R}_{>0}$. For $\phi$ in \eqref{eq:phi_class} with the corresponding $\psi$, we achieve $r<\rsoo{\phi}<\rho_{max}$ by prescribing a temporal behavior to $\rsntt{\psi}$ through the design parameters $\gamma$ and $\rho_{max}$ as
\begin{align}
	&-\gamma(t)+\rho_{max}< \rsntt{\psi}< \rho_{max}. \label{eq:inequality}
\end{align}

Note the use of $\rsntt{\psi}$ and not $\rsoo{\phi}$ itself. The connection between the non-temporal $\rsntt{\psi}$ and the temporal $\rsoo{\phi}$ is made by the performance function $\gamma$. In fact, $\gamma$ prescribes temporal behavior that, in combination with $\rsntt{\psi}$, mimics $\rsoo{\phi}$ as illustrated next.
\begin{example} 
	Fig. \ref{fig:principle_F} visualizes the idea for the eventually-operator $\phi_1:=F_{[0,\infty)}\psi_1$, while Fig. \ref{fig:principle_G} expresses the always-operator $\phi_2:=G_{[0,\infty)}\psi_2$. Note that these figures show the funnel in \eqref{eq:inequality}, hence imposing prescribed temporal behavior on $\rsntt{\psi}$. It is easy to verify that if $\rsntt{\psi_1}\in (-\gamma_1(t)+\rho_{1,max},\rho_{1,max})$ and $\rsntt{\psi_2}\in (-\gamma_2(t)+\rho_{2,max},\rho_{2,max})$ for all $t\in\mathbb{R}_{\ge 0}$ as in Fig. \ref{fig:11}, i.e. \eqref{eq:inequality} is satisfied,  then $\phi_1$ and $\phi_2$ are satisfied. For instance, in Fig. \ref{fig:principle_F} the lower funnel $-\gamma_1(t)+\rho_{1,max}$ forces $\rsntt{\psi_1} > r:=0$ by no later than approximately $4.5$ time units. Thus, the formulas $\phi_1:=F_{[0,\infty)}\psi_1$ or also $\phi_3:=F_{[2,5]}\psi_1$  are satisfied, which means that $\rsoo{\phi_1} > 0$ and $\rsoo{\phi_3} > 0$.
	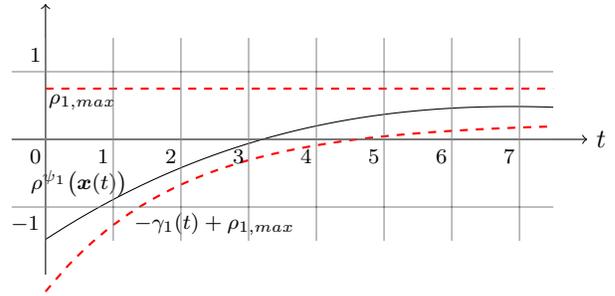
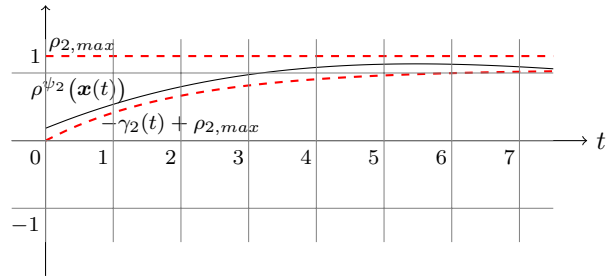
\begin{figure}
		\centering
		\begin{subfigure}[b]{0.48\textwidth}
			\begin{tikzpicture}[scale=0.9]
			\draw[->] (-0.5,0) -- (8,0) node[right] {$t$};
			\draw[->] (0,-2) -- (0,2) node[above] {};
			\draw[scale=1,domain=-0:7.5,smooth,variable=\x,red,dashed,thick] plot ({\x},{-((3-0.5)*exp(-0.5*\x)+0.5)+0.75});
			\draw[scale=1,domain=-0:7.5,smooth,variable=\x,red,dashed,thick] plot (\x,0.75);
			\draw[scale=1,domain=-0:7.5,smooth,variable=\x,black] plot (\x,{0.00177816*\x^3 - 0.0659376*\x^2 + 0.654718*\x - 1.47947});
			\node at (0.55,0.55) {\footnotesize $\rho_{1,max}$};
			\node at (2.5,-1.25) {\footnotesize $-\gamma_1(t)+\rho_{1,max}$};
			\node at (0.5,-0.65) {\footnotesize $\rsntt{\psi_1}$};
			\node at (-0.15,-0.25) {\footnotesize $0$};
			\node at (0.85,-0.25) {\footnotesize $1$};
			\node at (1.85,-0.25) {\footnotesize $2$};
			\node at (2.85,-0.25) {\footnotesize $3$};
			\node at (3.85,-0.25) {\footnotesize $4$};
			\node at (4.85,-0.25) {\footnotesize $5$};
			\node at (5.85,-0.25) {\footnotesize $6$};
			\node at (6.85,-0.25) {\footnotesize $7$};
			\node at (-0.15,1.25) {\footnotesize $1$};
			\node at (-0.3,-1.25) {\footnotesize $-1$};
			\draw[step=1,gray,very thin] (-0.5,-1.5) grid (7.5,1.5);
			\end{tikzpicture}\caption{Funnel for $\phi_1=F_{[0,\infty)}\psi_1$ s.t. $\rs{\phi_1}> r$ with $r:=0$}\label{fig:principle_F}
		\end{subfigure}
		\begin{subfigure}[b]{0.48\textwidth}
			\begin{tikzpicture}[scale=0.9]
			\draw[->] (-0.5,0) -- (8,0) node[right] {$t$};
			\draw[->] (0,-2) -- (0,2) node[above] {};
			\draw[scale=1,domain=-0:7.5,smooth,variable=\x,red,dashed,thick] plot ({\x},{-((1.25-0.2)*exp(-0.5*\x)+0.2)+1.25});
			\draw[scale=1,domain=-0:7.5,smooth,variable=\x,red,dashed,thick] plot (\x,1.25);
			\draw[scale=1,domain=-0:7.5,smooth,variable=\x,black] plot ({\x},{0.00188389*\x^3 - 0.0523356*\x^2 + 0.404143*\x + 0.183103});
			\node at (0.55,1.4) {\footnotesize $\rho_{2,max}$};
			\node at (2,0.25) {\footnotesize $-\gamma_2(t)+\rho_{2,max}$};
			\node at (0.5,0.75) {\footnotesize $\rsntt{\psi_2}$};
			\node at (-0.15,-0.25) {\footnotesize $0$};
			\node at (0.85,-0.25) {\footnotesize $1$};
			\node at (1.85,-0.25) {\footnotesize $2$};
			\node at (2.85,-0.25) {\footnotesize $3$};
			\node at (3.85,-0.25) {\footnotesize $4$};
			\node at (4.85,-0.25) {\footnotesize $5$};
			\node at (5.85,-0.25) {\footnotesize $6$};
			\node at (6.85,-0.25) {\footnotesize $7$};
			\node at (-0.15,1.25) {\footnotesize $1$};
			\node at (-0.3,-1.25) {\footnotesize $-1$};
			\draw[step=1,gray,very thin] (-0.5,-1.5) grid (7.5,1.5);
			\end{tikzpicture}\caption{Funnel for $\phi_2=G_{[0,\infty)}\psi_2$ s.t. $\rs{\phi_2}> r$ with $r:=0$}\label{fig:principle_G}
		\end{subfigure}
		\caption{Connection between $\rsnt{\psi}$ and $\rs{\phi}$}
		\label{fig:11}
	\end{figure}
\end{example}

The choice of the design parameters $\gamma$, $\rho_{max}$, and $r$ will be discussed in Section \ref{sec:strategy}. Therefore, define the global optimum of $\rsnt{\psi}$ as 
\begin{align}
\rso{\psi}:=\sup_{\boldsymbol{x}\in\mathbb{R}^n} \rsnt{\psi}.
\end{align} 

The function $\rsnt{\psi}$ is continuous and concave due to Assumption \ref{assumption:2} and \ref{assumption:4}, which makes the calculation of $\rso{\psi}$ straightforward. If $\rso{\psi}> 0$, it holds that $\phi$ is feasible, i.e., $\exists\boldsymbol{x}:\mathbb{R}_{\ge 0}\to \mathbb{R}^n$ s.t. $(\boldsymbol{x},0)\models \phi$.
\begin{assumption}\label{assumption:3}
	The optimum of $\rsnt{\psi}$ is s.t. $\rso{\psi}> 0$. 
\end{assumption}

Equation \eqref{eq:inequality} can now be written as
\begin{align}
	-\gamma(t) <\rsntt{\psi}-\rho_{max}< 0, \label{eq:inequality1}
\end{align}
which resembles \eqref{eq:constrained_funnel} by defining $M:=0$ and the one-dimensional error 
\begin{align}
	e(\boldsymbol{x}):=\rsnt{\psi}-\rho_{max}.\label{eq:e11}
\end{align}

Furthermore, define the normalized and the transformed error as
\begin{align}
\xi(\boldsymbol{x},t)&:=\frac{e(\boldsymbol{x})}{\gamma(t)},\label{eq:e12}\\
	\epsilon(\boldsymbol{x},t)&:=S\big(\xi(\boldsymbol{x},t)\big)=\ln\Big(-\frac{\xi(\boldsymbol{x},t)+1}{\xi(\boldsymbol{x},t)}\Big).\label{eq:e13}
\end{align}  

Hence, we can write \eqref{eq:inequality1} as $-\gamma(t)< e(t)< 0$, which in turn leads to $-1< \xi(t)< 0$. Applying the transformation function $S$ to this inequality finally results in $-\infty<\epsilon(t)<\infty$. In order to have a feasible problem, the condition $\xi\big(\boldsymbol{x}(0),0\big)\in\Omega_\xi:=(-1,0)$ needs to hold. As a notational rule, when talking about the solution $\boldsymbol{x}(t)$ of \eqref{system_noise} at time $t$, we use $e(t)$, $\xi(t)$, and $\epsilon(t)$, while we use $e(\boldsymbol{x})$, $\xi(\boldsymbol{x},t)$, and $\epsilon(\boldsymbol{x},t)$ when we talk about $\boldsymbol{x}$ as a state.

The second objective in this paper is to consider formulas $\theta$ as in \eqref{eq:theta_class}, called sequential formulas. The name stems from the fact, that the atomic temporal formulas contained in \eqref{eq:theta1_class} and \eqref{eq:theta2_class} can be processed sequentially. Therefore, the hybrid system framework of \cite{goebel2012hybrid} will be used. We are now ready for the formal problem definition:

\begin{problem}\label{problem}
	Consider the system given in \eqref{system_noise} subject to a STL formula $\theta$ as in \eqref{eq:theta_class}. Design a piecewise-continuous feedback control law $\boldsymbol{u}(\boldsymbol{x},t)$ such that $r<\rsoo{\theta}<\rho_{max}$.  
\end{problem}

Note that $\theta$ boils down to an atomic temporal formula $\phi$ as in \eqref{eq:phi_class} if $N=1$, i.e., $\theta$ is a superset of $\phi$. Our problem solution consists of a three-step procedure: First, a continuous feedback control law $\boldsymbol{u}(\boldsymbol{x},t)$ is designed in Theorem \ref{theorem:3} such that \eqref{eq:inequality} is satisfied, which means that $\rsnt{\psi}$ follows a prescribed behavior. Second, $\gamma$ is designed in Theorem \ref{lemma:gamma} such that $r<\rsoo{\phi}<\rho_{max}$ if $\boldsymbol{u}(\boldsymbol{x},t)$ from Theorem \ref{theorem:3} is used. Third, Theorem \ref{theorem:hybrid_system} states a hybrid control strategy such that $r<\rsoo{\theta}<\rho_{max}$. Section \ref{sec:strategy} covers Theorem \ref{theorem:3} and \ref{lemma:gamma} and hence achieves satisfaction of atomic temporal formulas, i.e., $r<\rsoo{\phi}<\rho_{max}$, while Section \ref{sec:sequencing} covers Theorem \ref{theorem:hybrid_system} and leads to satisfaction of sequential formulas, i.e., $r<\rsoo{\theta}<\rho_{max}$.

\section{Control Law for Atomic Temporal Formulas}
\label{sec:strategy}

As explained previously, in a first step we derive a control law $\boldsymbol{u}(\boldsymbol{x},t)$ such that $\rsntt{\psi}$ satisfies \eqref{eq:inequality}, while in a second step $\gamma$ is designed such that $\rsoo{\phi}> r$. Recall \eqref{eq:e11}, \eqref{eq:e12}, and \eqref{eq:e13}, then the dynamics of $\epsilon$ are given by $\dot{\epsilon} = \frac{\partial \epsilon}{\partial \xi} \dot{\xi} =-\frac{1}{\gamma\xi(1+\xi)}(\frac{\partial \rsnt{\psi}}{\partial\boldsymbol{x}}^T\dot{\boldsymbol{x}}-\xi\dot{\gamma})$ since $\frac{\partial \epsilon}{\partial \xi} = -\frac{1}{\xi(1+\xi)}$ and $\dot{\xi}=\frac{1}{\gamma}(\dot{e}-\xi\dot{\gamma})$. Note that $\dot{e} = \frac{\partial e(\boldsymbol{x})}{\partial\boldsymbol{x}}^T\dot{\boldsymbol{x}}$ with $\frac{\partial e(\boldsymbol{x})}{\partial\boldsymbol{x}}=\frac{\partial \rsnt{\psi}}{\partial\boldsymbol{x}}$. 

\begin{theorem}\label{theorem:3}
Consider the system \eqref{system_noise} and a formula $\phi$ as in \eqref{eq:phi_class} with the corresponding $\psi$. If $\xi\big(\boldsymbol{x}_0,0\big)\in \Omega_\xi:=(-1,0)$, $\rho_{max}\in\big(\max\big(0,\rsnto{\psi}\big),\rso{\psi}\big)$, and Assumptions \ref{assumption:1}-\ref{assumption:3} are satisfied, then the control law 
\begin{align}\label{equ:control}
\boldsymbol{u}(\boldsymbol{x},t):=-\epsilon(\boldsymbol{x},t)g^T(\boldsymbol{x})\frac{\partial \rsnt{\psi}}{\partial \boldsymbol{x}}
\end{align}
guarantees that \eqref{eq:inequality} is satisfied for all $t\in\mathbb{R}_{\ge 0}$ with all closed-loop signals being well-posed, i.e., continuous and bounded.

\begin{proof}
We proceed as follows: in the first step (Step A), we apply Lemma \ref{theorem:sontag1} and show that there exists a maximal solution $\xi(t)$ such that $\xi(t)\in\Omega_\xi$ for all $t\in[0,\tau_{max})=\mathcal{J}\subseteq \mathbb{R}_{\ge 0}$. The second step (step B) consists of using Lemma \ref{theorem:sontag2} to show that $\tau_{max}=\infty$, which proves the main result. 


Step A: First, define the stacked vector $\boldsymbol{y}:=\begin{bmatrix}
\boldsymbol{x}^T & \xi
\end{bmatrix}^T$. Consider the closed-loop system that is obtained by inserting \eqref{equ:control} into \eqref{system_noise} resulting in $\dot{\boldsymbol{x}}:=H_1(\boldsymbol{x},\xi)=f(\boldsymbol{x})-\ln(-\frac{\xi+1}{\xi})g(\boldsymbol{x})g^T(\boldsymbol{x})\frac{\partial \rsnt{\psi}}{\partial\boldsymbol{x}}+\boldsymbol{w}$. We also obtain $\dot{\xi}:=H_2(\boldsymbol{x},\xi,t)=\frac{1}{\gamma(t)}\big(\frac{\partial \rsnt{\psi}}{\partial\boldsymbol{x}} H_1(\boldsymbol{x},\xi)-\xi \dot{\gamma}(t)  \big)$, which results in $\dot{\boldsymbol{y}}:=H(\boldsymbol{y},t)=\begin{bmatrix}
H_1(\boldsymbol{x},\xi) & H_2(\boldsymbol{x},\xi,t)
\end{bmatrix}^T$. According to the assumptions, it holds that $\boldsymbol{x}_0$ is such that $\xi(\boldsymbol{x}_0,0)\in\Omega_\xi=(-1,0)$, which is non-empty and open. Next, define the time-varying and non-empty set $\Omega_{\boldsymbol{x}}(t):=\{\boldsymbol{x}\in\mathbb{R}^n|-1<\xi(\boldsymbol{x},t)=\frac{\rsnt{\psi}-\rho_{max}}{\gamma(t)}<0\}$, which has the property that for $t_1<t_2$ it is true that $\Omega_{\boldsymbol{x}}(t_2)\subseteq \Omega_{\boldsymbol{x}}(t_1)$ since $\gamma(t)$ is non-increasing in $t$. Note that $\Omega_{\boldsymbol{x}}(t)$ is bounded due to Assumption \ref{assumption:4}. We denote $\Omega_{\boldsymbol{x}_0} :=\Omega_{\boldsymbol{x}}(0)$ and remark that $\boldsymbol{x}_0\in\Omega_{\boldsymbol{x}}(0)$. Due to \cite[Proposition 1.4.4]{aubin2009set}, the following holds: if a function is continuous, then the inverse image of an open set under this function is open. With $\xi_0(\boldsymbol{x})=\xi(\boldsymbol{x},0)$, it holds that the inverse image ${\xi_0}^{-1}(\Omega_\xi)=\Omega_{\boldsymbol{x}_0}$ is open. Note therefore, that $\rsnt{\psi}$ is a continuously differentiable function due to Assumption \ref{assumption:2}. Finally, define the open, bounded, and non-empty set $\Omega_{\boldsymbol{y}}:=\Omega_{\boldsymbol{x}_0} \times \Omega_\xi$, which does not depend on $t$. It consequently holds that $\boldsymbol{y}_0=\begin{bmatrix}
\boldsymbol{x}_0^T & \xi_0
\end{bmatrix}^T\in\Omega_{\boldsymbol{y}}$.

Next, the conditions in Lemma \ref{theorem:sontag1} for the IVP $\dot{\boldsymbol{y}}=H(\boldsymbol{y},t)$ with $\boldsymbol{y}_0\in\Omega_{\boldsymbol{y}}$ and $H(\boldsymbol{y},t):\Omega_{\boldsymbol{y}}\times \mathbb{R}_{\ge 0} \to \mathbb{R}^{n+1}$ need to be checked: 1) $H(\boldsymbol{y},t)$ is locally Lipschitz on $\boldsymbol{y}$ since $f(\boldsymbol{x})$, $g(\boldsymbol{x})$, and $\epsilon=\ln\big(-\frac{\xi+1}{\xi}\big)$ are  locally Lipschitz continuous on $\boldsymbol{y}$ for each $t\in \mathbb{R}_{\ge 0}$. This also holds for $\frac{\partial \rsnt{\psi}}{\partial\boldsymbol{x}}$ due to Assumption \ref{assumption:2}. 2) $H(\boldsymbol{y},t)$ is continuous on $t$ for each fixed $\boldsymbol{y}\in\Omega_{\boldsymbol{y}}$ due to continuity of $\gamma(t)$ and $\dot{\gamma}(t)$. Finally, $\Omega_{\boldsymbol{y}}$ is non-empty and open. Applying Lemma \ref{theorem:sontag1}, there exists a maximal solution with $\boldsymbol{y}(t)\in\Omega_{\boldsymbol{y}}$ for all $t\in[0,\tau_{max})=\mathcal{J}\subseteq \mathbb{R}_{\ge 0}$ and $\tau_{max}>0$. Consequently, there exist $\xi(t)\in\Omega_\xi$ and $\boldsymbol{x}(t)\in\Omega_{\boldsymbol{x}_0}$ for all $t \in \mathcal{J}$. 

Step B: From Step A), it is known that $\boldsymbol{y}(t)\in\Omega_{\boldsymbol{y}}$ for all $t\in[0,\tau_{max})=\mathcal{J}$. Next, we show that $\tau_{max}=\infty$ by contradiction of Lemma \ref{theorem:sontag2}. Therefore, assume $\tau_{max}<\infty$ and consider the Lyapunov function $V(\epsilon)=\frac{1}{2}\epsilon^2$.  Hence, it holds that 
\begin{align}\label{eq:p1}
\dot{V}=\epsilon\dot{\epsilon} =\epsilon\Big(-\frac{1}{\gamma\xi(1+\xi)}\big(\frac{\partial \rsnt{\psi}}{\partial \boldsymbol{x}}^T\dot{\boldsymbol{x}}-\xi\dot{\gamma}\big)\Big).
\end{align} 
Inserting \eqref{system_noise} into \eqref{eq:p1} results in 
\begin{align}\label{eq_p3}
\dot{V}=-\frac{\epsilon}{\gamma\xi(1+\xi)}\Big(\frac{\partial \rsnt{\psi}}{\partial \boldsymbol{x}}^T\big(f(\boldsymbol{x})+g(\boldsymbol{x})\boldsymbol{u}+\re{\boldsymbol{w}}\big)-\xi\dot{\gamma}\Big).
\end{align} 
Define $\alpha(t)=-\frac{1}{\gamma\xi(1+\xi)}$ which satisfies $\alpha(t)\in[\frac{4}{\gamma_0},\infty)\in\mathbb{R}_{> 0}$ for all $t\in \mathcal{J}$. This follows since $\frac{4}{\gamma_0}\le -\frac{1}{\gamma_0\xi(1+\xi)}\le-\frac{1}{\gamma\xi(1+\xi)}\le-\frac{1}{\gamma_\infty\xi(1+\xi)}< \infty$ for $\xi\in\Omega_\xi$. Next, \eqref{eq_p3} can be upper bounded as 
\begin{small}
\begin{align}
\dot{V}&\le |\epsilon| \alpha \big(\|\frac{\partial \rsnt{\psi}}{\partial \boldsymbol{x}}\|\|f(\boldsymbol{x})+\re{\boldsymbol{w}}\|+|\xi\dot{\gamma}|\big)+\epsilon \alpha\frac{\partial \rsnt{\psi}}{\partial \boldsymbol{x}}^Tg(\boldsymbol{x})\boldsymbol{u}\\
&\le |\epsilon| \alpha k_1 +\epsilon \alpha\frac{\partial \rsnt{\psi}}{\partial \boldsymbol{x}}^Tg(\boldsymbol{x})\boldsymbol{u},\label{eq:p2}
\end{align}
\end{small}where the last inequality and the positive constant $k_1$ derives as follows: it holds that $\|f(\boldsymbol{x})\|<\infty$ and $\|\frac{\partial \rsnt{\psi}}{\partial \boldsymbol{x}}\|<\infty$ since $\boldsymbol{x}(t)\in\Omega_{\boldsymbol{x}_0}$ for all $t\in \mathcal{J}$ and due to the extreme value theorem and continuity of $f(\boldsymbol{x})$ and $\frac{\partial \rsnt{\psi}}{\partial \boldsymbol{x}}$. Furthermore, $\boldsymbol{w}$ and $\dot{\gamma}$ are bounded. Next, insert the control law \eqref{equ:control} into \eqref{eq:p2}, which results in 
\begin{align}
\dot{V}\le |\epsilon| \alpha k_1 -\epsilon^2 \alpha\frac{\partial \rsnt{\psi}}{\partial \boldsymbol{x}}^Tg(\boldsymbol{x})g^T(\boldsymbol{x})\frac{\partial \rsnt{\psi}}{\partial \boldsymbol{x}}
\end{align} 
\begin{align}
\le|\epsilon|\alpha \big(k_1-|\epsilon| \lambda_{min}(g(\boldsymbol{x})g^T(\boldsymbol{x}))\|\frac{\partial \rsnt{\psi}}{\partial \boldsymbol{x}}\|^2\big),
\end{align} 
where $\lambda_{min}(g(\boldsymbol{x})g^T(\boldsymbol{x}))>0$ is the minimum eigenvalue of $g(\boldsymbol{x})g^T(\boldsymbol{x})$, which is positive according to Assumption \ref{assumption:1}. It holds that $\|\frac{\partial \rsnt{\psi}}{\partial \boldsymbol{x}}\|^2\ge k_2>0$ for a positive constant $k_2$ since $\rsnt{\psi}$ is concave as a result of Assumption \ref{assumption:4}, and hence $\frac{\partial \rsnt{\psi}}{\partial \boldsymbol{x}}=0$ if and only if $\rsnt{\psi}=\rso{\psi}$. However, this case has been excluded since $\rsnt{\psi} \in (-\gamma(t)+\rho_{max},\rho_{max})$ for all $t\in \mathcal{J}$, which ensures that $\rsnt{\psi}<\rso{\psi}$ due to the assumption that $\rho_{max}\in\big(\max\big(0,\rsnto{\psi}\big),\rso{\psi}\big)$. Finally, $\dot{V}$ can be upper bound as 
\begin{align}
\dot{V}\le|\epsilon|\alpha \big(k_1-|\epsilon| \lambda_{min}(g(\boldsymbol{x})g^T(\boldsymbol{x}))k_2\big).
\end{align}
Hence, $\dot{V}\le 0$ if $\frac{k_1}{\lambda_{min}(g(\boldsymbol{x})g^T(\boldsymbol{x}))k_2}\le|\epsilon|$ and it can be concluded that the transformed error $|\epsilon|$ will be upper bounded due to the level sets of $V(\epsilon)$ as $|\epsilon(t)|\le \max\left(|\epsilon(0)|,\frac{k_1}{\lambda_{min}(g(\boldsymbol{x})g^T(\boldsymbol{x}))k_2}\right)$, which leads to the conclusion that $\epsilon(t)$ is upper and lower bounded by some constants $\epsilon_u$ and $\epsilon_l$, respectively. In other words, it holds that $\epsilon_l\le \epsilon(t)\le \epsilon_u$. By using the inverse of $S(\cdot)$, the normalized error $\xi(t)$ can be bounded by $-1<\xi_l:=-\frac{1}{\exp(\epsilon_l+1)}\le \xi(t) \le \xi_u:=-\frac{1}{\exp(\epsilon_u+1)}<0$, which means that $\xi(t)\in [\xi_l,\xi_u]=:\Omega_{\xi}^\prime\subset \Omega_\xi$ for all $t\in \mathcal{J}$. Recall \eqref{eq:e12} and note that if $\xi(t)$ evolves in a compact set, then $\rsntt{\psi}$ will evolve in a compact set $\Omega_\rho^\prime:=[\rho_l,\rho_u]$ for some constants $\rho_l$ and $\rho_u$. Again, due to \cite[Proposition 1.4.4]{aubin2009set} it holds that the inverse image ${\rho^{\psi}}^{-1}(\Omega_\rho^\prime)=\{\boldsymbol{x}\in \Omega_{\boldsymbol{x}}|\rho_l\le \rsnt{\psi}\le \rho_u\}=:\Omega_{\boldsymbol{x}}^\prime$ is closed and also bounded due to Assumption \ref{assumption:4}, which hence excludes finite escape time of the state $\boldsymbol{x}$. Consequently, it can be concluded that $\boldsymbol{x}(t)$ evolves in a compact set, i.e., $\boldsymbol{x}(t)\in \Omega_{\boldsymbol{x}}^\prime\subset\Omega_{\boldsymbol{x}_0}$ for all $t\in \mathcal{J}$. Define the compact set $\Omega_{\boldsymbol{y}}^\prime:=\Omega_{\boldsymbol{x}}^\prime\times \Omega_{\xi}^\prime$ and notice that $\Omega_{\boldsymbol{y}}^\prime\subset \Omega_{\boldsymbol{y}}$ by which it follows that there is no $t\in \mathcal{J}=[0,\tau_{max})$ such that $\boldsymbol{y}\notin \Omega_{\boldsymbol{y}}^\prime$. By contradiction of Lemma \ref{theorem:sontag2} it follows that $\tau_{max}=\infty$, i.e., $\mathcal{J}=\mathbb{R}_{\ge 0}$. 

The control law $\boldsymbol{u}(\boldsymbol{x},t)$ is well-posed, i.e., continuous and bounded, because $\rsnt{\psi}$ is approximated by a smooth function, while $\epsilon(\boldsymbol{x},t)$ and $g(\boldsymbol{x})$ are  locally Lipschitz continuous on $\boldsymbol{x}$. Due to the extreme value theorem, these functions are bounded on $\boldsymbol{x}$. Also, $\gamma(t)$ is continuous with $0<\gamma(t)<\infty$. It follows that all closed-loop signals are well-posed. 
\end{proof}
\end{theorem}

The second step is to show that the control law \eqref{equ:control} in Theorem \ref{theorem:3} results in $r<\rsoo{\phi}<\rho_{max}$ if $\gamma$ is properly designed. The variable $t_*\in\mathbb{R}_{\ge 0}$  is s.t.
\begin{align}
t_*\in\begin{cases} a &\text{ if } \phi=G_{[a,b]}\psi \\ 
[a,b] &\text{ if } \phi=F_{[a,b]}\psi,
\end{cases}
\end{align} 
which will enforce $r<\rsntt{\psi}<\rho_{max}$ for all $t\ge t_*$ by the choice of $\gamma$ in the remainder. This consequently leads to  $r<\rsoo{\phi}<\rho_{max}$ by the choice of $t_*$. We select $r\in[0,\rho_{max})$ and define feasibility of a formula $\phi$ with respect to $r$, $\boldsymbol{x}_0$, and $t_*$. \ref{def:form_feasibility}. 
\begin{definition}\label{def:form_feasibility}
	A formula $\phi$ as in \eqref{eq:phi_class} is feasible with respect to $r$, $\boldsymbol{x}_0$, and $t_*$ if and only if: 1) $t_*>0$ or 2) $t_*=0$ and $\rsnto{\psi}> r$.
\end{definition}

For the design of $\gamma$ assume that $\phi$ is feasible w.r.t. $r$, $\boldsymbol{x}_0$, and $t_*$ and recall that $\gamma(t):=(\gamma_0-\gamma_\infty)\exp({-lt})+\gamma_\infty$. The crucial part of Theorem \ref{theorem:3} is the assumption that $\xi(\boldsymbol{x}_0,0)\in\Omega_\xi$. It is possible to choose $\gamma_0$ such that $\xi(\boldsymbol{x}_0,0)\in\Omega_\xi$, which is equivalent to $-1<\frac{\rsnto{\psi}-\rho_{max}}{\gamma(0)}<0$. It should also hold that $-\gamma_0+\rho_{max}\ge r$ if $t_*=0$ due to \eqref{eq:inequality} and since we want $r<\rsntt{\psi}$ for all $t\ge t_*$. This is illustrated in Fig. \ref{fig:principle_G} with $t_*=0$ (since $\phi_2=G_{[0,\infty)}\psi_2$) and $r:=0$ and where it should hence hold that $-\gamma_0+\rho_{max}\ge 0$ is satisfied as indicated by the dashed line. To conclude, $\gamma_0$ is 
\begin{align}\label{eq:g1}
\gamma_{0}\in\begin{cases}(\rho_{max}-\rsnto{\psi},\infty) &\text{if } t_*>0\\
(\rho_{max}-\rsnto{\psi},\rho_{max}-r] &\text{if } t_*=0.\end{cases}
\end{align}

At $t=\infty$, it is required that  $\max(-\gamma_0+\rho_{max},r)\le -\gamma_\infty+\rho_{max}<\rho_{max}$, where the left inequality enforces that $-\gamma+\rho_{max}$ is a non-decreasing function, which in turn leads to $\gamma$ being non-increasing. The right inequality stems from \eqref{eq:inequality}. Therefore, we set 
\begin{align}\label{eq:g2}
\gamma_{\infty}\in \Big(0,\min\big(\gamma_0,\rho_{max}-r\big)\Big].
\end{align}

The smaller $\gamma_{\infty}$ is selected, the tighter the funnel will be as $t\to\infty$. For the calculation of $l$, three cases need to be distinguished: 1) $\rsnto{\psi}> r$, 2) $\rsnto{\psi}\le r$ and $t_*>0$, and 3) $\rsnto{\psi}\le r$ and $t_*=0$. Case 3) can be excluded since $\phi$ is feasible w.r.t. $r$, $\boldsymbol{x}_0$, and $t_*$. Next, select $l$ as
\begin{align}\label{eq:g3}
l\in \begin{cases}
\mathbb{R}_{\ge 0} & \text{ if } -\gamma_0+\rho_{max}\ge r\\
-\frac{\ln\big(\frac{r+\gamma_{\infty}-\rho_{max}}{-(\gamma_{0}-\gamma_{\infty})}\big)}{t_*} & \text{ if } -\gamma_0+\rho_{max}< r, t_*>0,
\end{cases}
\end{align} 
which ensures that $-\gamma(t_*)+\rho_{max}\ge r$. Under \eqref{equ:control}, this consequently leads to $\rsntt{\psi}>r$ for all $t\ge t_*$ since $\gamma$ is non-increasing. 

\begin{theorem}\label{lemma:gamma}
Consider the system \eqref{system_noise} and a formula $\phi$ as in \eqref{eq:phi_class}. If Assumptions \ref{assumption:1}-\ref{assumption:3} hold, $r\in[0,\rho_{max})$, the control law in \eqref{equ:control} is used, and  $\phi$ is feasible w.r.t. $r$, $\boldsymbol{x}_0$, and $t_*$, then choosing $\gamma_0$, $\gamma_\infty$, and $l$ as in \eqref{eq:g1}, \eqref{eq:g2}, and \eqref{eq:g3}, respectively, ensures that $0\le r<\rsoo{\phi}<\rho_{max}$, i.e., $(\boldsymbol{x},0)\models \phi$.

\begin{proof}
Choosing $\gamma_0$ as in \eqref{eq:g1} ensures $\xi(\boldsymbol{x}_0,0)\in\Omega_\xi$, while additionally choosing $\gamma_\infty$ and $l$ as in \eqref{eq:g2} and \eqref{eq:g3} ensures $\rsoo{\phi}> r$ if \eqref{equ:control} is applied. This follows since by the above choice, we impose $-\gamma(t_*)+\rho_{max}= r$ for case 2) while case 1) already has $\rsnto{\psi}> r$. Note for case 2) that solving the equation $-\gamma(t_*)+\rho_{max}=:r$ for $l$ results in $l=-\frac{\ln\big(\frac{r+\gamma_{\infty}-\rho_{max}}{-(\gamma_{0}-\gamma_{\infty})}\big)}{t_*}$. Hence, the control law \eqref{equ:control} enforces $\rsnts{\psi}> r$, which consequently leads to $\rsoo{\phi}>r$ due to the choice of $t_*$. It hence holds that $r<\rsoo{\phi}<\rho_{max}$.
\end{proof}
\end{theorem}
\begin{remark}\label{remark:2}
The assumption of feasibility w.r.t. $r$, $\boldsymbol{x}_0$, and $t_*$ is a necessary assumption. However, if a formula is not feasible w.r.t $r$, $\boldsymbol{x}_0$, and $t_*$, the formula can be relaxed as discussed in \cite{ghosh2016diagnosis}. 
\end{remark}
\begin{remark}
In combination, Theorem \ref{theorem:3} and \ref{lemma:gamma} provide a control strategy such that $(\boldsymbol{x},0) \models \phi$. However, a steep performance function $\gamma$ might result in a high control effort. Therefore, it may in practice be useful to choose $t_*$ as big and $l$ as small as possible.
\end{remark}
\section{Control strategy for sequential formulas}
\label{sec:sequencing}
In this section, we develop a hybrid control strategy for sequential formulas $\theta$ as in \eqref{eq:theta_class}, which either correspond to $\theta^{s_1}$ or $\theta^{s_2}$ as in \eqref{eq:theta1_class} or \eqref{eq:theta2_class}, respectively. Note that both of these consist of $N$ atomic temporal formulas: $\theta^{s_1}$ entails $N$ atomic temporal formulas $\phi_i$ with $[a_i,b_i]$ for all  $i\in\{1,\hdots,N\}$. Similarly, $\theta^{s_2}$ boils down to $N-1$ atomic temporal formulas $\phi_i=F_{[a_i,b_i]}\psi_i$ with $i\in\{1,\hdots,N-1\}$, $a_{i}:=\sum_{k=1}^ic_k$, $b_i:=\sum_{k=1}^id_k$,  and $\phi_N$. For instance, $F_{[c_1,d_1]}\big(\psi_1 \wedge F_{[c_2,d_2]}(\psi_2 \wedge F_{[c_3,d_3]}\psi_3)\big)$ is satisfied if and only if $F_{[c_1,d_1]}\psi_1 \wedge F_{[c_1+c_2,d_1+d_2]}\psi_2 \wedge  F_{[c_1+c_2+c_3,d_1+d_2+d_3]}\psi_3:=F_{[a_1,b_1]}\psi_1 \wedge F_{[a_2,b_2]}\psi_2 \wedge  F_{[a_3,b_3]}\psi_3$ is satisfied. To conclude, $\theta$ consists of $N$ atomic temporal formulas $\phi_i$ with $i\in\{1,\hdots,N\}$. Each $\phi_i$ entails a robustness function denoted by $\rsnt{\psi_i}$ and corresponding design parameters $t_{i,*}$, $r_i$, $\rho_{i,max}$, and $\gamma_i(t)=(\gamma_{i,0}-\gamma_{i,\infty})\exp(-l_{i}t)+\gamma_{i,\infty}$ in accordance with $t_*$, $r$, $\rho_{max}$, and $\gamma$ in Section \ref{sec:strategy}.  Each $\phi_i$ will be processed one at a time. If $\phi_i$ has been satisfied, the next atomic temporal formula $\phi_{i+1}$ becomes active and a switch takes place. Denote the time sequence of these switching times by $\{\Delta_1:=0,\Delta_2,\hdots,\Delta_N\}$ where $\Delta_i\le \Delta_{i+1}$. Note that $t_{i,*}$, $r_i$, $\rho_{i,max}$, $\gamma_{i,0}$, $\gamma_{i,\infty}$, and $l_i$ need to be calculated during runtime at each switching time $\Delta_i$. Furthermore, set $p:=\begin{cases} 1  \text{ if } \theta=\theta^{s_1}\\ 0 \text{ if } \theta=\theta^{s_2} \end{cases}$  and $m_i:=\begin{cases} 1 \text{ if } \phi_i=G_{[a_i,b_i]}\psi_i \\ 0  \text{ if } \phi_i=F_{[a_i,b_i]}\psi_i. \end{cases}$

A hybrid control strategy in the framework introduced in Definition \ref{def:hybrid_system} will be used to process each $\phi^i$ sequentially.
\begin{definition}\cite{goebel2012hybrid}\label{def:hybrid_system}
	A hybrid system is a tuple $\mathcal{H}:=(C,F,D,G)$, where $C$, $D$, $F$, and $G$ are the flow and jump set and the possibly set-valued flow and jump map, respectively. The discrete and continuous dynamics are 
	\begin{align*}
	\begin{cases}
	\dot{\boldsymbol{z}}\in F(\boldsymbol{z}) & \text{if }\boldsymbol{z} \in C  \\
	\boldsymbol{z}^+\in G(\boldsymbol{z}) & \text{if }\boldsymbol{z} \in D. 
	\end{cases}
	\end{align*}
\end{definition}

Define $\boldsymbol{p}_f:=\begin{bmatrix} t_* & r & \rho_{max} & \gamma_{0} & \gamma_{\infty} & l \end{bmatrix}^T$, gathering all parameters defining the funnel in \eqref{eq:inequality}, and the hybrid state $\boldsymbol{z}:=\begin{bmatrix} q & \boldsymbol{x}^T & t & \Delta & \boldsymbol{p}_f^T  \end{bmatrix}^T\in \{1,\hdots,N+1\}\times \mathbb{R}^n \times \mathbb{R}_{\ge 0}^8 =:\mathcal{Z}$. Note that $\Delta$ is the value of the latest switching time. In adherence to the terminology in \cite{goebel2012hybrid}, we interchangeably call switches jumps. The discrete state $q$ indicates which formula $\phi_q$ is currently active, while $q=N+1$ indicates the final discrete state when $\theta$ has already been satisfied. In the proof of Theorem \ref{theorem:3}, it was shown that $\boldsymbol{x}(t)\in\Omega_{\boldsymbol{x}}^\prime$ for all $t\in\mathbb{R}_{\ge 0}$, where $\Omega_{\boldsymbol{x}}^\prime$ is a compact set.  Let $\Omega_{q,\boldsymbol{x}}^\prime$ denote $\Omega_{\boldsymbol{x}}^\prime$ corresponding to the formula $\phi_q$. Next, define the sets $\mathcal{X}_q:=\{\boldsymbol{x}\in\mathbb{R}^n|r_q<\rsnt{\psi_q}<\rho_{q,max}\}$ and $\mathcal{Y}_q:= [0,{b_q}^p+{\big(\sum_{i=1}^qd_i\big)}^{1-p}-1] \times t_{q,*} \times r_q \times \rho_{q,max} \times \gamma_{q,0} \times \gamma_{q,\infty} \times l_{q}$. Note that $p$ determines if $[a_q,b_q]$ or $[c_q,d_q]$ is used. For all $q\in\{1,\hdots,N\}$, set
\begin{small}
	\begin{align}\label{eq:t_star}
	t_{q,*}\in
	\begin{cases}
	a_q &\text{if } p=1, m_q=1\\
	[a_q,b_q] &\text{if } p=1, m_q=0\\
	c_q &\text{if } p=0, m_q=1\\
	[c_q,d_q] &\text{if } p=0, m_q=0.\\
	\end{cases}
	\end{align} 
\end{small}
Next, define the set $\mathcal{D}_{q}$ that indicates satisfaction of $\phi_q$ and leads to a jump to process $\phi_{q+1}$. For $q\in\{1,\hdots,N\}$, define 
\begin{small}
	\begin{align*}
	\mathcal{D}_{q}:= \begin{cases}
	q \times \mathcal{X}_q \times ({b_q}^p+{d_q}^{1-p}-1-p\Delta)  \times \mathcal{Y}_q \hspace{0.69cm} \text{if } m_q=1\\ 
	q \times \mathcal{X}_q \times ([{a_q}^p+{c_q}^{1-p}-1,t_{q,*}]-p\Delta) \times  \mathcal{Y}_q \;\text{if } m_q=0,
	\end{cases}
	\end{align*} 
\end{small}which indicates that $\rsq{\phi_q}>r_q$ if $\boldsymbol{z}\in\mathcal{D}_{q}$. This follows since $\boldsymbol{x}\in\mathcal{X}_q$ at $t={b_q}^p+{d_q}^{1-p}-1-p\Delta$ for $m_q=1$ or $\boldsymbol{x}\in\mathcal{X}_q$ at $t\in ([{a_q}^p+{c_q}^{1-p}-1,t_{q,*}]-p\Delta)$ for $m_q=0$ under the control law \eqref{equ:control} indicates that $\phi_q$ is satisfied.  Note that $\Delta$ only takes effect if $p=1$ ($\theta=\theta^{s_1}$) to ensure that $\phi_q$ is satisfied within $[a_q,b_q]$, while for $p=0$ ($\theta=\theta^{s_2}$) the formula $\phi_{q+1}$ is directly processed next when $\phi_q$ is satisfied. Further define $\mathcal{D}_{N+1}:=(N+1) \times \Omega_{N,\boldsymbol{x}}^\prime \times T  \times \mathcal{Y}_N$ for $T:={b_N}^p+(\sum_{i=1}^Nd_i)^{1-p}-1$, which is needed for a technical reason in the proof of Theorem \ref{theorem:hybrid_system}. Similarly, define the continuous domain $\mathcal{C}_q$ for $q\in\{1,\hdots,N\}$ as 
\begin{small}
	\begin{align*}
	\mathcal{C}_{q}:=
	\begin{cases}
	q\times \Omega_{q,\boldsymbol{x}}^\prime \times [0, {b_q}^p+{d_q}^{1-p}-1-p\Delta] \times \mathcal{Y}_q \;\;\text{if } m_q=1\\
	q \times \text{cl}(\Omega_{q,\boldsymbol{x}}^\prime\setminus \mathcal{X}) \times [0, t_{q,*}-p\Delta] \times \mathcal{Y}_q \hspace{1.05cm}\text{if } m_q=0,
	\end{cases}
	\end{align*}
\end{small}where cl($\cdot$) denotes the closure. Also define $\mathcal{C}_{N+1}:=(N+1)\times \Omega_{N,\boldsymbol{x}}^\prime\times [0,T]  \times \mathcal{Y}_N$.  Finally, the jump and flow sets are given by   
\begin{align*}
D&:=  \cup_{i=1}^{N+1}  \mathcal{D}_{i}  \\
C&:= \cup_{i=1}^{N+1} \mathcal{C}_{i}.
\end{align*} 
The flow map is given by 
\begin{align*}
F:=\begin{bmatrix} 0 & {\big(f(\boldsymbol{x})+g(\boldsymbol{x})\boldsymbol{u}_q+\boldsymbol{w}\big)}^T & 1 & {\boldsymbol{0}_7}^T\end{bmatrix}^T
\end{align*} 
with the control law in \eqref{equ:control} as $\boldsymbol{u}_{q}=-\epsilon_q g^T(\boldsymbol{x}) \frac{\partial \rsnt{\psi_q}}{\partial \boldsymbol{x}}$ for all $q\in \{1,\hdots,N\}$ and $\boldsymbol{u}_{N+1}=-\epsilon_N g^T(\boldsymbol{x}) \frac{\partial \rsnt{\psi_N}}{\partial \boldsymbol{x}}$, where $\epsilon_q$  corresponds to $\epsilon$  based on $\phi_q$. By abbreviating $q^\prime:=q+1$, define $\boldsymbol{p}_s(q):=\begin{bmatrix}
q^\prime & \boldsymbol{x}^T & 0 & \Delta_{q^\prime} & t_{q^\prime,*} & r_{q^\prime} & \rho_{q^\prime,max} &\gamma_{q^\prime,0} & \gamma_{q^\prime,\infty} &l_{q^\prime} 
\end{bmatrix}^T$. Then, the jump map is given by
\begin{small}
	\begin{align*}
	G:=
	\begin{cases}
	{\boldsymbol{p}_s(q)} \hspace{3.22cm}
	\text{if}\; q\notin\{N,N+1\}, \boldsymbol{z}\in D \\ 
	\begin{bmatrix}
	N+1 & \boldsymbol{x}^T & \boldsymbol{0}_2^T & {\boldsymbol{p}_f}^T 
	\end{bmatrix}^T 
	\;\text{if}\; q\in\{N,N+1\}, \boldsymbol{z}\in D,  
	\end{cases}
	\end{align*}
\end{small}where we set $\Delta_{q^\prime}:=\Delta+t$, accumulating the elapsed time after jumps. Select $t_{q^\prime,*}$ as in \eqref{eq:t_star},  $\rho_{q^\prime,max}\in \big(\max(0,\rsnt{\psi_{q^\prime}}),\rso{\psi_{q^\prime}}\big)$, and $r_{q^\prime}\in[0,\rho_{q^\prime,max})$ as in the assumptions of Theorem \ref{theorem:3} and \ref{lemma:gamma}, respectively. The parameters $\gamma_{i,0}$, $\gamma_{i,\infty}$, and $l_i$ need to be chosen as in \eqref{eq:g1}, \eqref{eq:g2}, and \eqref{eq:g3}: $\gamma_{q^\prime,0}\in\begin{cases}(\rho_{q^\prime,max}-\rsnt{\psi_{q^\prime}},\infty)&\text{if } t_{q^\prime,*}-p\Delta_{q^\prime}>0\\ (\rho_{q^\prime,max}-\rsnt{\psi_{q^\prime}},\rho_{q^\prime,max}-r_{q^\prime}]&\text{if } t_{q^\prime,*}-p\Delta_{q^\prime}=0 \end{cases}$, $\gamma_{q^\prime,\infty}\in \Big(0,\min\big(\gamma_{q^\prime,0},\rho_{q^\prime,max}-r_{q^\prime}\big)\Big]$, and 
\begin{small}
\begin{align*}
l_{q^\prime}\in
\begin{cases}
\mathbb{R}_{\ge 0} \hspace{3.025cm} \text{if }  -\gamma_{q^\prime,0}+\rho_{q^\prime,max}\ge r_{q^\prime}\\
-\frac{\ln\big(\frac{r_{q^\prime}+\gamma_{q^\prime,\infty}-\rho_{q^\prime,max}}{-(\gamma_{q^\prime,0}-\gamma_{q^\prime,\infty})}\big)}{t_{q^\prime,*}-p\Delta_{q^\prime}} \;\;\; \text{if: }\bullet -\gamma_{q^\prime,0}+\rho_{q^\prime,max}< r_{q^\prime},\\ \hspace{4.06cm} \bullet \;t_{q^\prime,*}-p\Delta_{q^\prime}>0.
\end{cases}
\end{align*} 
\end{small}Note that $G$ is hence a set-valued map. The initial state is set to $\boldsymbol{z}_0:=\begin{bmatrix}
1 & {\boldsymbol{x}_0}^T & 0 & 0 & t_{1,*} & r_{1} & \rho_{1,max} &\gamma_{1,0} & \gamma_{1,\infty} & l_{1}
\end{bmatrix}^T$.  Now, we are ready to state the main result of this section.

\begin{theorem}\label{theorem:hybrid_system}
	
	Consider the system \eqref{system_noise} and a formula $\theta$ as in \eqref{eq:theta_class}. The hybrid system $\mathcal{H}=(C,F,D,G)$ results in $r:=\min(r_1,\hdots,r_N)<\rsoo{\theta}<\rho_{max}:=\min(\rho_{1,max},\hdots,\rho_{N,max})$, i.e., $(\boldsymbol{x},0)\models\theta$, if each $\phi_q$ in $\theta$ is feasible w.r.t. $r_q$, $\boldsymbol{x}(\Delta_q)$, and $t_{q,*}+|p-1|\Delta_q$.

\begin{proof} 
First, note that the third case in the proof of Theorem \ref{lemma:gamma} is again excluded by the assumption of feasibility w.r.t. $r_q$, $\boldsymbol{x}(\Delta_q)$, and $t_{q,*}+|p-1|\Delta_q$. To show that $\theta$ is satisfied, we need to show that eventually $\phi_N$ is satisfied. Therefore, we show that the compact set $\mathcal{A}=(N+1)\times \bigcup \limits_{i=1}^{N}\Omega_{i,\boldsymbol{x}}^\prime \times [0,T]^3 \times [0,\max_i r_i]  \times [0,\max_i \rho_{i,max}] \times [0,\max_i \gamma_{i,0}] \times [0,\max_i \gamma_{i,\infty}] \times [0,\max_i l_{i}]$ is asymptotically stable. A hybrid Lyapunov-function candidate is $V(\boldsymbol{z}):=\big(q-(N+1)\big)^2$, which is positive on $(C\cup D) \setminus \mathcal{A}$. During flows it is easy to see that $\dot{V}=0$, while during jumps $V(\boldsymbol{z}^+)-V(\boldsymbol{z})=\big(q+1-(N+1)\big)^2-\big(q-(N+1)\big)^2=\big(q-N\big)^2-\big(q-(N+1)\big)^2<0$ for $q\in\{1,\hdots,N\}$. According to the invariance principle in \cite[Theorem 23]{goebel2009hybrid}, we now need to show that no complete solution can stay in $V(\boldsymbol{z})=\mu>0$. This is true due to the following fact: for each state $q=\{1,\hdots,N\}$, the control law  $\boldsymbol{u}_{q}(t)=-\epsilon_q g^T(\boldsymbol{x}) \frac{\partial \rsnt{\psi_q}}{\partial \boldsymbol{x}}$ of Theorem \ref{theorem:3} is applied to the system. Furthermore, $r_{q}$ and $\rho_{q,max}$ are chosen as in Theorem \ref{theorem:3} and \ref{lemma:gamma}, while $\gamma_{q,0}$, $\gamma_{q,\infty}$, and $l_{q}$ are chosen accordingly. This guarantees that each $\phi_q$ is satisfied with $\rho^{\phi_q}(\boldsymbol{x},\Delta_q)> r_q$. Hence, $\phi_q$ will eventually be satisfied and lead to a jump due to the structure of the jump set $D$, which decreases $V(\boldsymbol{z})$. Note that $\Delta_{q}$ ensures that each formula $\phi_q$ in $\theta^{s_1}$ is satisfied within $[a_q,b_q]$, while each $\phi_q$ in $\theta^{s_2}$ is processed without the use of $\Delta_q$. Hence, note that each solution is complete and it can be concluded that $\mathcal{A}$ is asymptotically stable, which leads to the conclusion that $\theta$ is satisfied with $r:=min(r_1,\hdots,r_N)<\rs{\theta}<\rho_{max}:=\min(\rho_{1,max},\hdots,\rho_{N,max})$.
\end{proof}
\end{theorem}

\section{Simulations}
\label{sec:simulations}

\begin{figure}[b]
\centering
\input{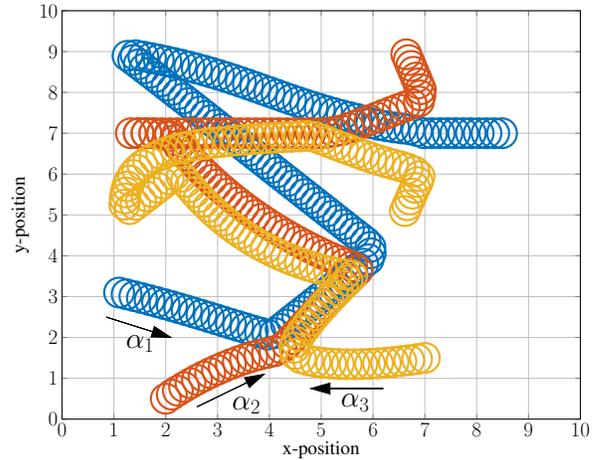}\caption{Continuous trajectory for $\phi_1$, $\phi_2$, $\phi_3$, $\phi_4$}\label{fig:5}
\end{figure}

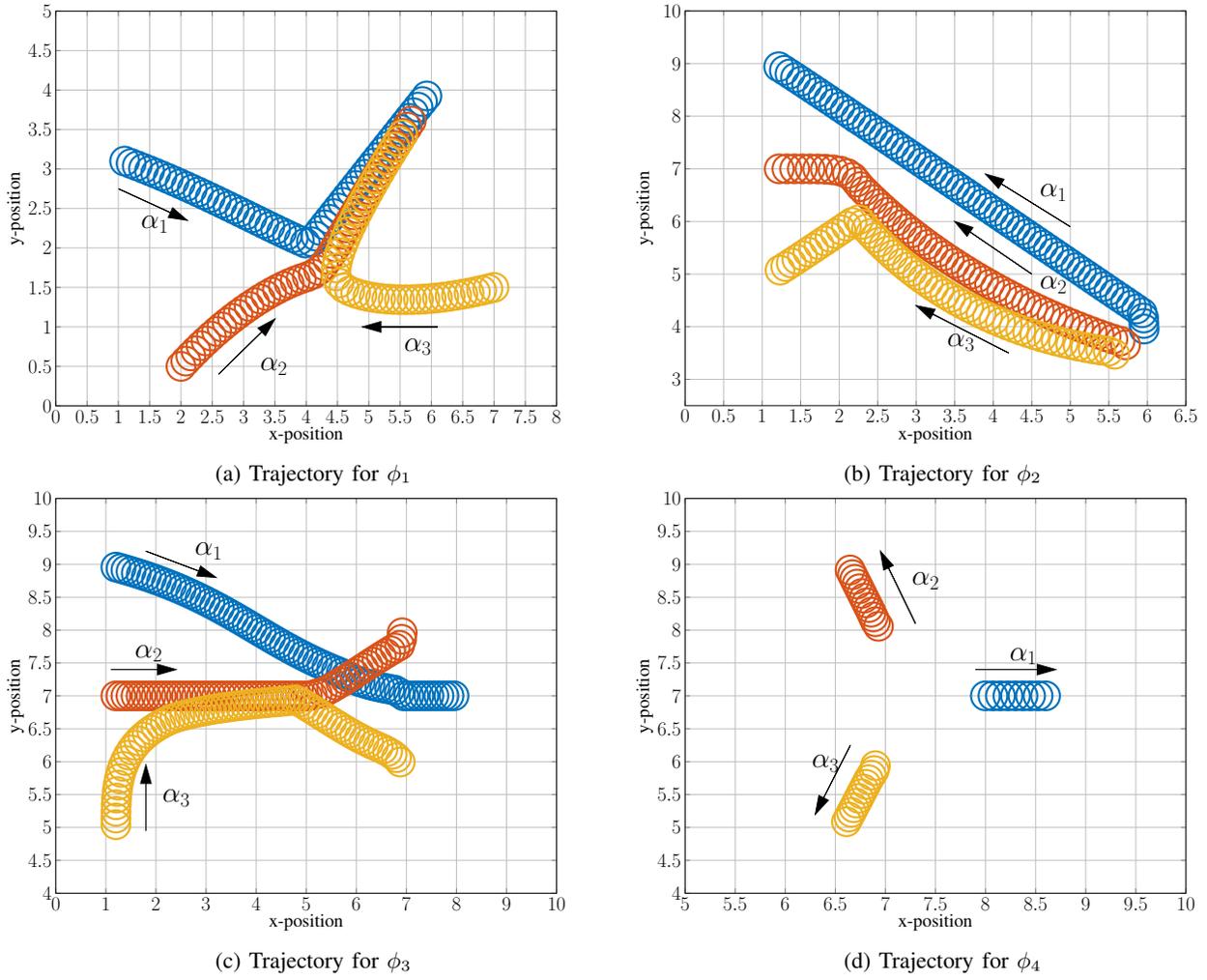
\begin{figure*}
\centering
\begin{subfigure}{0.48\textwidth}
%
%
\definecolor{mycolor1}{rgb}{0.00000,0.44700,0.74100}%
\definecolor{mycolor2}{rgb}{0.85000,0.32500,0.09800}%
\definecolor{mycolor3}{rgb}{0.92900,0.69400,0.12500}%
\begin{tikzpicture}[scale=0.45]

\begin{axis}[%
width=6.028in,
height=4.754in,
at={(1.011in,0.642in)},
scale only axis,
xmin=0,
xmax=8,
xlabel={x-position},
xmajorgrids,
ymin=0,
ymax=5,
ylabel={y-position},
ymajorgrids,
axis background/.style={fill=white},
title style={font=\LARGE},xlabel style={font=\LARGE},ylabel style={font=\LARGE},legend style={font=\LARGE},ticklabel style={font=\LARGE},
]
\addplot [color=mycolor1,line width=1.8pt,mark size=12.5pt,only marks,mark=o,mark options={solid},forget plot]
  table[row sep=crcr]{%
1.1	3.1\\
1.20349113288927	3.0692054843373\\
1.3015962191778	3.03927731075184\\
1.40242187852948	3.00780030972585\\
1.49812462107609	2.97727287304403\\
1.59600130357357	2.94542121387616\\
1.69544964102894	2.91242889464723\\
1.79591813939566	2.87847908484225\\
1.89111931529444	2.8457581621786\\
1.98701126814562	2.81228067969177\\
2.08315317023827	2.77821632016415\\
2.17914778442614	2.74372950464544\\
2.2746403199806	2.70897803981766\\
2.36931702380351	2.67411193655393\\
2.46702972967345	2.63772824534669\\
2.56293327168161	2.60165603500035\\
2.65686494823527	2.56601189840789\\
2.7521154294217	2.52958843355095\\
2.84794744462499	2.49270245088191\\
2.94371578331155	2.45564959584856\\
3.0388688155834	2.41870074734303\\
3.13294742396092	2.38209941344319\\
3.22799556722768	2.34512169754954\\
3.32321559411832	2.30816023344331\\
3.41795367697951	2.27156215471031\\
3.51169626382615	2.23562592987932\\
3.60583812692219	2.19993115040539\\
3.69971988322183	2.16485757049218\\
3.79436890342028	2.13016733581862\\
3.88910319655955	2.09627000984875\\
3.98472027058642	2.0630509443361\\
4.0839094874953	2.07869680739607\\
4.15468722521663	2.14981888553987\\
4.22586938432377	2.22122920102038\\
4.29702797790912	2.29255316584207\\
4.36772080232884	2.36337995004386\\
4.43878366395774	2.43456058142028\\
4.50998423192929	2.50586919250753\\
4.58132660987077	2.57731389794765\\
4.65218033961703	2.64826582905225\\
4.72347854652075	2.71966094509505\\
4.79465564510089	2.79093382600896\\
4.86586398051978	2.86223757862551\\
4.93656803342109	2.93303625539626\\
5.00744266983953	3.00400574085786\\
5.07838123272536	3.07503925075432\\
5.14917615112964	3.14592891176645\\
5.21998594967807	3.21683344561585\\
5.29072524469889	3.2876673909772\\
5.36146186681323	3.35849879494582\\
5.43224154774703	3.42937370383307\\
5.50292680577806	3.50015487751906\\
5.5738381443309	3.57116393423422\\
5.64461645105569	3.64204232617019\\
5.71543253387452	3.71296279872939\\
5.78616371148201	3.78380542497568\\
5.85687466932739	3.85464079501073\\
5.92754028192903	3.92545872741844\\
};
\addplot [color=mycolor2,line width=1.8pt,mark size=12.5pt,only marks,mark=o,mark options={solid},forget plot]
  table[row sep=crcr]{%
2	0.5\\
2.07713401901991	0.566584222219921\\
2.15648915412793	0.632685473879285\\
2.23766198829031	0.697946879652795\\
2.32024764998184	0.762040195153916\\
2.40026806188997	0.8220339744044\\
2.48460801288656	0.883113893039006\\
2.56918291192305	0.942224660619049\\
2.65362119288003	0.999166934410441\\
2.73756922050383	1.05377484668974\\
2.82359044927933	1.10769765879756\\
2.9109371569638	1.1603733905655\\
2.99887890144053	1.21130985481437\\
3.08671659253615	1.26008970857676\\
3.17606951201536	1.3075537911852\\
3.26586309152041	1.35303822672954\\
3.35704394543836	1.39691548182665\\
3.44822516576177	1.43840676412949\\
3.54139067846473	1.47825174179982\\
3.63450629052161	1.51538589711761\\
3.72839256221445	1.54995179731169\\
3.82409994997171	1.58201047231377\\
3.92043901950631	1.61077363808398\\
4.01735116882382	1.63615442538078\\
4.11233357685645	1.66804243467822\\
4.20001236868594	1.71614068350484\\
4.2773353200084	1.77994358201651\\
4.34583176826556	1.85335908252975\\
4.40882938310683	1.9312774708014\\
4.46901020006439	2.01134092559906\\
4.52792347691201	2.09249224732447\\
4.58613682343582	2.17382343405061\\
4.6442786596985	2.25529153490925\\
4.70250466957548	2.33661264716441\\
4.7612303622908	2.41809218391261\\
4.8201617479164	2.49917941843126\\
4.87930881082307	2.57982813807133\\
4.93895338480757	2.66040242520195\\
4.99901050551684	2.74078888548784\\
5.05935441203805	2.82083808821889\\
5.12008517448713	2.9007085772248\\
5.18095378895248	2.98010311398505\\
5.24214532443855	3.05929371611711\\
5.30368372502315	3.13833668110803\\
5.36542697344661	3.21707376737817\\
5.42744839669588	3.29562020391548\\
5.48970071583037	3.37393372075111\\
5.55224526320668	3.45210492355022\\
5.61500205314944	3.53004315924512\\
5.67797400697358	3.60775732405732\\
};
\addplot [color=mycolor3,line width=1.8pt,mark size=12.5pt,only marks,mark=o,mark options={solid},forget plot]
  table[row sep=crcr]{%
7	1.5\\
6.90171361871839	1.48085667335974\\
6.80231784294921	1.4626125699773\\
6.70223461314498	1.44539973315861\\
6.60187446922837	1.42933361854203\\
6.50163361274127	1.41451247843293\\
6.401891228796	1.40101704146343\\
6.29929963695692	1.38848178624161\\
6.19827959828365	1.37752027229903\\
6.09581583714102	1.36785895463874\\
5.99588484415947	1.3599115866226\\
5.89578375284839	1.35348129047305\\
5.79352927541157	1.34857826001503\\
5.69295848753322	1.34549461583622\\
5.59234951687072	1.34424826993548\\
5.49066740315509	1.34499421319682\\
5.38951385868175	1.34790318010795\\
5.28857357997341	1.35316534850324\\
5.18812277111276	1.36099853155433\\
5.08758129246668	1.37176803561144\\
4.98757983091005	1.38584115760452\\
4.88881431193236	1.40367425992025\\
4.7911819842147	1.42613500454331\\
4.69481146684918	1.45503932624682\\
4.60294912875788	1.49479023584221\\
4.52370120784916	1.55615953869622\\
4.48056796945944	1.64645538640708\\
4.47998143731025	1.74668737941237\\
4.50269913869827	1.84421797494055\\
4.53715379058045	1.93838513942429\\
4.57841817285434	2.02994082406525\\
4.62407374406371	2.11942113661593\\
4.67279239943318	2.2072306577899\\
4.72380838827713	2.29372304119547\\
4.77655462842583	2.37903324221912\\
4.83072232130903	2.46340124814194\\
4.88616205900056	2.54710678576045\\
4.9425501808812	2.63003986072841\\
4.9998690497119	2.71246291898631\\
5.05785366511776	2.79421975021013\\
5.11648416824132	2.87546392692743\\
5.17569613029455	2.95624976619151\\
5.23536143051124	3.0365213327725\\
5.29569791351885	3.11666474520092\\
5.35634019417078	3.19627289512209\\
5.41738025703123	3.27553600304454\\
5.47881386447611	3.35450323295347\\
5.54060869286645	3.43317761632816\\
};

\addplot[area legend,solid,draw=black,fill=black,forget plot]
table[row sep=crcr] {%
x	y\\
1.41597240929956	2.59873730570925\\
1	2.75\\
1	2.75\\
1	2.75\\
1	2.75\\
1	2.75\\
1.83194481859911	2.4474746114185\\
1.79667322329366	2.38656532772068\\
2.1	2.35\\
1.86721641390457	2.50838389511633\\
1.83194481859911	2.4474746114185\\
}--cycle;
\node[right, align=left, text=black]
at (axis cs:1.3,2.3) {\Huge $\alpha_1$};

\addplot[area legend,solid,draw=black,fill=black,forget plot]
table[row sep=crcr] {%
x	y\\
2.94475752061372	0.668144738255113\\
2.6	0.4\\
2.6	0.4\\
2.6	0.4\\
2.6	0.4\\
2.6	0.4\\
3.28951504122743	0.936289476510226\\
3.23027570446329	0.984117270776321\\
3.5	1.1\\
3.34875437799158	0.88846168224413\\
3.28951504122743	0.936289476510226\\
}--cycle;
\node[right, align=left, text=black]
at (axis cs:3.2,0.5) {\Huge $\alpha_2$};

\addplot[area legend,solid,draw=black,fill=black,forget plot]
table[row sep=crcr] {%
x	y\\
5.64746543778802	1\\
6.1	1\\
6.1	1\\
6.1	1\\
6.1	1\\
6.1	1\\
5.19493087557604	1\\
5.19493087557604	0.932983842066303\\
4.9	1\\
5.19493087557604	1.0670161579337\\
5.19493087557604	1\\
}--cycle;
\node[right, align=left, text=black]
at (axis cs:5.5,0.75) {\Huge $\alpha_3$};
\end{axis}
\end{tikzpicture}%
\caption{Trajectory for $\phi_1$}\label{fig:1}
\end{subfigure}
\begin{subfigure}{0.48\textwidth}
%
%
\definecolor{mycolor1}{rgb}{0.00000,0.44700,0.74100}%
\definecolor{mycolor2}{rgb}{0.85000,0.32500,0.09800}%
\definecolor{mycolor3}{rgb}{0.92900,0.69400,0.12500}%
\begin{tikzpicture}[scale=0.45]

\begin{axis}[%
width=6.028in,
height=4.754in,
at={(1.011in,0.642in)},
scale only axis,
xmin=0,
xmax=6.5,
xlabel={x-position},
xmajorgrids,
ymin=2.5,
ymax=10,
ylabel={y-position},
ymajorgrids,
axis background/.style={fill=white},
title style={font=\LARGE},xlabel style={font=\LARGE},ylabel style={font=\LARGE},legend style={font=\LARGE},ticklabel style={font=\LARGE},
]
\addplot [color=mycolor1,line width=1.8pt,mark size=12.5pt,only marks,mark=o,mark options={solid},forget plot]
  table[row sep=crcr]{%
5.95798491753768	3.95598989769185\\
5.95463186537511	4.06108619215902\\
5.95111017570152	4.16508363212664\\
5.93166736232734	4.26476847576699\\
5.86047697932767	4.33709359763759\\
5.78951974251282	4.4080126953341\\
5.71735697685483	4.47974995391814\\
5.64578402440987	4.55046586810984\\
5.57378794639819	4.62132890017844\\
5.50186060880597	4.69287629467855\\
5.4307119921999	4.76383452739797\\
5.35882898872548	4.83552880330141\\
5.28753834579641	4.9066353099192\\
5.21561120564668	4.97837933846005\\
5.14435363409699	5.04945778730261\\
5.07256023783772	5.12107267897996\\
5.0015132275554	5.1919446772639\\
4.93003223098831	5.26325095062694\\
4.85818229743808	5.33492636190735\\
4.78720509957115	5.40573199336663\\
4.71596148286368	5.47680399998473\\
4.64451758008494	5.54807618988185\\
4.57294003959266	5.61948186185974\\
4.50129577852962	5.69095405213492\\
4.42965173590289	5.76242578113182\\
4.3580746279809	5.83383029790871\\
4.28663070827816	5.90510131995487\\
4.21538553419393	5.97617326630138\\
4.14440374214004	6.04698148211975\\
4.07274938909866	6.11845942726227\\
4.0015240582733	6.18950803816832\\
3.92983058137152	6.26102206515845\\
3.85872876316173	6.33194416766946\\
3.78736048618055	6.40313010790282\\
3.71584629059953	6.47445943085805\\
3.6443057908048	6.54581260857012\\
3.57285713134409	6.61707158448332\\
3.50161645254184	6.68812030817113\\
3.43069737083565	6.75884525536391\\
3.35943022159892	6.82991399424828\\
3.28799303790361	6.9011486703679\\
3.21656106691888	6.97237421943518\\
3.14530573570026	7.04341940256884\\
3.07439364506481	7.11411781358921\\
3.00332746622941	7.18496493120879\\
2.93233614928399	7.25573210403689\\
2.86103391991362	7.32680342666159\\
2.78971022731532	7.39788987318689\\
2.71864435100061	7.46871261278667\\
2.64757331625711	7.53953320262448\\
2.57632358067592	7.6105239415964\\
2.50528593851777	7.68129480799435\\
2.43437924399377	7.75192599830777\\
2.36317851481395	7.82284004026714\\
2.29223391620924	7.89348813804215\\
2.21898589619579	7.96186615810795\\
2.14887614957807	8.03355770687582\\
2.07811999606557	8.10437666066063\\
2.00707377103368	8.174851189307\\
1.9358286647671	8.24518187465622\\
1.8645393471738	8.31534069035587\\
1.79300386158141	8.38558122561118\\
1.72143682203688	8.45570685726369\\
1.64978931105103	8.52575113438404\\
1.57816306015232	8.59557551717932\\
1.50638627524402	8.66527859158722\\
1.43438433190514	8.73481404283655\\
1.36212073201864	8.80399818425905\\
1.28934136519618	8.87260205013457\\
1.21542737459024	8.93996427177015\\
};
\addplot [color=mycolor2,line width=1.8pt,mark size=12.5pt,only marks,mark=o,mark options={solid},forget plot]
  table[row sep=crcr]{%
5.71627103808696	3.65478145550188\\
5.6454413813499	3.72716395374875\\
5.55265211537372	3.76511330669811\\
5.46298642802168	3.80949237999484\\
5.37290299363945	3.85575565340347\\
5.28373645423689	3.90316930142287\\
5.19554033028792	3.95163960558284\\
5.1083652265921	4.00107834415043\\
5.02103642240909	4.05212704640214\\
4.93485270047362	4.10400371620821\\
4.84866776317048	4.15737607617865\\
4.76374369831355	4.21144122276503\\
4.67896416158594	4.26688828155537\\
4.59555071598066	4.32289702399603\\
4.51241847216326	4.38017379398382\\
4.43074533165474	4.43788470997531\\
4.34947928709855	4.49674960788511\\
4.26869957384729	4.55671543220451\\
4.18848333453911	4.61772716996263\\
4.108905336404	4.6797279901442\\
4.03003771656714	4.74265939587974\\
3.95194975574781	4.80646138711614\\
3.87470768023668	4.87107263157036\\
3.79837449159574	4.93643064190169\\
3.72300982316695	5.00247195719624\\
3.64779646955717	5.06992459815742\\
3.5737068359419	5.13792505848152\\
3.50078985839219	5.20640599295535\\
3.42828760789292	5.27608042313506\\
3.35631178624985	5.34686498991145\\
3.28572550109257	5.41790345896661\\
3.21582917272225	5.48988778786716\\
3.14672136897258	5.56272636324616\\
3.07849569655665	5.63632484850398\\
3.01124045819779	5.71058659127295\\
2.94440109024028	5.78614188759831\\
2.87873839286282	5.86214191155809\\
2.81373201616935	5.93919198254347\\
2.75008193764427	6.01645833614911\\
2.68729987422206	6.09452106479113\\
2.62550272651174	6.17323872354643\\
2.56429964510235	6.25312698104653\\
2.50433267101461	6.33335452905694\\
2.44523530058353	6.41440571520235\\
2.38715037083387	6.4960898812862\\
2.32979697499319	6.57881540513823\\
2.27374563536305	6.66175731727785\\
2.22538752277172	6.74980774218471\\
2.17751570941169	6.83768978133132\\
2.10708788407642	6.9087577084991\\
2.0185488423002	6.9554703686435\\
1.92169456591835	6.98071701157953\\
1.82225833200859	6.99262048161965\\
1.72226101667017	6.99757767081664\\
1.62223356869977	6.9993549861501\\
1.52219425318962	6.99987388017078\\
1.42213406741304	6.99998504356129\\
1.32212381495506	6.99999930873196\\
1.22210699291244	6.99999999778848\\
};
\addplot [color=mycolor3,line width=1.8pt,mark size=12.5pt,only marks,mark=o,mark options={solid},forget plot]
  table[row sep=crcr]{%
5.57543157643603	3.47719663403341\\
5.49277586790875	3.53358531008846\\
5.39469819711168	3.55450889423305\\
5.29641633627867	3.57792000192765\\
5.19929957284845	3.60349353259358\\
5.10217877630825	3.63153347527335\\
5.00638854568454	3.66164675192081\\
4.91078818924221	3.69417981579324\\
4.81667146487714	3.72867686306555\\
4.72292749523443	3.7655253948149\\
4.62966292368992	3.80470696859181\\
4.53810171411814	3.84567990924971\\
4.44718169294508	3.88888791420419\\
4.35699925329112	3.93429618805774\\
4.26869957384729	3.98128862552129\\
4.1812725254461	4.03035909648772\\
4.09480168260857	4.08145885505794\\
4.00936640610853	4.13453402571733\\
3.92504162096695	4.18952584072211\\
3.84189764023788	4.24637089477014\\
3.76000003160006	4.30500141586785\\
3.67940952342914	4.36534555104439\\
3.60018194685001	4.42732766535165\\
3.52236821023724	4.49086865242582\\
3.44601430271875	4.55588625476865\\
3.37116132342047	4.62229539183313\\
3.29708611864463	4.69072407458907\\
3.22462839084142	4.76037642584845\\
3.15381420278278	4.83115775083536\\
3.08397846896296	4.90369858714377\\
3.01587184771301	4.97717836579047\\
2.9488666718781	5.05222504879796\\
2.88303957941124	5.12873627490499\\
2.81905299518321	5.20587694291674\\
2.7563330146583	5.2842654353909\\
2.69493796363099	5.36378660246868\\
2.63491969930548	5.44432116768616\\
2.57632358067592	5.5257463469883\\
2.519188487444	5.60793648824758\\
2.46308574979178	5.69146203639612\\
2.40854338113809	5.77547941927596\\
2.35515918797401	5.86053741398554\\
2.30300406663205	5.94647433751082\\
2.25222639900209	6.03262969208024\\
2.16619217001991	5.98156128247628\\
2.09495026486168	5.91113989530033\\
2.02359905060254	5.84083346942186\\
1.95221978574133	5.77062524281405\\
1.88084760813056	5.70051698844982\\
1.80945916187781	5.63047962119283\\
1.73797322816315	5.56044657737364\\
1.66638396731732	5.49044178465735\\
1.59474453629906	5.42056307739986\\
1.52295348856549	5.35077966447663\\
1.45102132990843	5.28120945616209\\
1.37883479805215	5.21193215444465\\
1.30620497560531	5.14315516666048\\
1.23260491363379	5.07541034928856\\
};

\addplot[area legend,solid,draw=black,fill=black,forget plot]
table[row sep=crcr] {%
x	y\\
4.55176745202975	6.30748413451841\\
5	5.9\\
5	5.9\\
5	5.9\\
5	5.9\\
5	5.9\\
4.1035349040595	6.71496826903682\\
4.06726788065338	6.62958631849517\\
3.9	6.9\\
4.13980192746562	6.80035021957847\\
4.1035349040595	6.71496826903682\\
}--cycle;
\node[right, align=left, text=black]
at (axis cs:4.55,6.55) {\Huge $\alpha_1$};

\addplot[area legend,solid,draw=black,fill=black,forget plot]
table[row sep=crcr] {%
x	y\\
4.09891525436963	5.40108474563037\\
4.5	5\\
4.5	5\\
4.5	5\\
4.5	5\\
4.5	5\\
3.69783050873927	5.80216949126073\\
3.65905486926628	5.71918050816293\\
3.5	6\\
3.73660614821226	5.88515847435853\\
3.69783050873927	5.80216949126073\\
}--cycle;
\node[right, align=left, text=black]
at (axis cs:4.55,4.8) {\Huge $\alpha_2$};

\addplot[area legend,solid,draw=black,fill=black,forget plot]
table[row sep=crcr] {%
x	y\\
3.70662096613339	3.87003427539996\\
4.2	3.5\\
4.2	3.5\\
4.2	3.5\\
4.2	3.5\\
4.2	3.5\\
3.21324193226678	4.24006855079992\\
3.18189467317062	4.1506145468376\\
3	4.4\\
3.24458919136294	4.32952255476224\\
3.21324193226678	4.24006855079992\\
}--cycle;
\node[right, align=left, text=black]
at (axis cs:3.35,3.7) {\Huge $\alpha_3$};
\end{axis}
\end{tikzpicture}%
\caption{Trajectory for $\phi_2$}\label{fig:2}
\end{subfigure}
\begin{subfigure}{0.48\textwidth}
%
%
\definecolor{mycolor1}{rgb}{0.00000,0.44700,0.74100}%
\definecolor{mycolor2}{rgb}{0.85000,0.32500,0.09800}%
\definecolor{mycolor3}{rgb}{0.92900,0.69400,0.12500}%
\begin{tikzpicture}[scale=0.45]

\begin{axis}[%
width=6.028in,
height=4.754in,
at={(1.011in,0.642in)},
scale only axis,
xmin=0,
xmax=10,
xlabel={x-position},
xmajorgrids,
ymin=4,
ymax=10,
ylabel={y-position},
ymajorgrids,
axis background/.style={fill=white},
title style={font=\LARGE},xlabel style={font=\LARGE},ylabel style={font=\LARGE},legend style={font=\LARGE},ticklabel style={font=\LARGE},
]
\addplot [color=mycolor1,line width=1.8pt,mark size=12.5pt,only marks,mark=o,mark options={solid},forget plot]
  table[row sep=crcr]{%
1.2007176764374	8.95749372787454\\
1.30368012320551	8.93606863668148\\
1.40159373535524	8.91487804663962\\
1.50321501119514	8.89202547170412\\
1.6077753706645	8.86757617878976\\
1.70712692828871	8.84344282661349\\
1.80883579283127	8.81780171116646\\
1.91228353538265	8.79072473005461\\
2.0105322553598	8.76405107285959\\
2.11011500487076	8.73603695715677\\
2.21054852763744	8.7067592535271\\
2.31140215641196	8.67629733362185\\
2.4071208663226	8.64637909236512\\
2.50308874130183	8.61537783145581\\
2.59899064651051	8.58337536288537\\
2.69455728913446	8.55045489264603\\
2.78956351278705	8.51670066368825\\
2.88382612024073	8.48219760099386\\
2.98098885820748	8.44558392986604\\
3.07677175628451	8.40846482546513\\
3.17113913385326	8.37092655117686\\
3.26408370470503	8.33305435476026\\
3.35869972176875	8.29363722690416\\
3.45165405211331	8.25413023976664\\
3.54582287450152	8.21339933110898\\
3.6382649226708	8.17281726583092\\
3.731655345821	8.13132641386802\\
3.82337829767029	8.09021208376812\\
3.91592088261557	8.04849178525376\\
4.00915655089883	8.00635231379943\\
4.10085107744259	7.96494054482755\\
4.19323639949001	7.92338334488834\\
4.28624060968862	7.88184954525132\\
4.37788685935602	7.84134175274023\\
4.47016376021223	7.80109010054408\\
4.56300139632986	7.761239410918\\
4.65632618756486	7.72192509799579\\
4.75005752788804	7.68327258039701\\
4.84410548291117	7.64539682086003\\
4.93836950723863	7.60840199171902\\
5.03273809061317	7.57238126339284\\
5.12708920302326	7.53741671153688\\
5.22129138427661	7.50357933714249\\
5.31658322949442	7.4704582634773\\
5.41267760321618	7.43819837756809\\
5.50800283972571	7.40733469759565\\
5.6036459243475	7.37751004400007\\
5.69996743442082	7.34742737392742\\
5.79183301486666	7.30627782204835\\
5.88580821015718	7.27081368201936\\
5.98123642091496	7.24017576853205\\
6.07799875464409	7.21321722558966\\
6.1755947570777	7.18928464584917\\
6.27375177590405	7.1678703415172\\
6.37231722640743	7.14858108914356\\
6.47117814875796	7.13111415059351\\
6.57020151230783	7.11523848678019\\
6.66919025070041	7.10077956486792\\
6.76839477847138	7.08693960895398\\
6.86023537628815	7.04591344725659\\
6.94948529902962	7.00059951028886\\
7.04971850382606	7.00065979145411\\
7.14978326906637	7.00084293591326\\
7.24985171719687	7.00096966379759\\
7.34985794860146	7.00106320873538\\
7.44996638657872	7.00113602784858\\
7.55015283001974	7.00119302362453\\
7.65031274392323	7.00123664782856\\
7.7504203143741	7.00127280911445\\
7.85049094821512	7.00130169551493\\
7.95052866359483	7.00131556549087\\
};
\addplot [color=mycolor2,line width=1.8pt,mark size=12.5pt,only marks,mark=o,mark options={solid},forget plot]
  table[row sep=crcr]{%
1.2007176764374	6.99999999999317\\
1.30076526132343	6.99999999999499\\
1.40079189942625	6.99999999999563\\
1.50135131924295	6.99999999999606\\
1.60177634307874	6.99999999999639\\
1.70274395800637	6.99999999999663\\
1.80388835215021	6.99999999999683\\
1.90413930031656	6.99999999999699\\
2.00428993458104	6.99999999999714\\
2.10536397919889	6.99999999999722\\
2.20563855918056	6.9999999999973\\
2.30636649610723	6.99999999999738\\
2.407415316811	6.99999999999747\\
2.50869596052057	6.99999999999752\\
2.6087849151846	6.99999999999752\\
2.70903442866365	6.99999999999752\\
2.80940986449408	6.99999999999752\\
2.90987971472696	6.99999999999752\\
3.01041190699856	6.99999999999752\\
3.1109714413241	6.99999999999752\\
3.2115190296986	6.99999999999752\\
3.31201047084206	6.99999999999752\\
3.41234943313946	6.99999999999752\\
3.51301088786821	6.99999999999752\\
3.61319379315976	6.99999999999752\\
3.71377087271964	6.99999999999752\\
3.81418116118865	6.99999999999752\\
3.91481768430551	6.99999999999752\\
4.01588335347395	6.99999999999752\\
4.11623885344829	6.99999999999752\\
4.21702226456403	6.99999999999752\\
4.31810920363519	6.99999999999752\\
4.41811621284527	6.99999999999752\\
4.51924014192076	6.99999999999752\\
4.62003894606489	6.99999999999764\\
4.72026387185728	7.00000000002487\\
4.82084126202842	7.00000000594414\\
4.92148169684698	7.00000125021247\\
5.02188622811669	7.00024164470012\\
5.12158849888157	7.01051472840567\\
5.2196781958686	7.03188987890964\\
5.31797541227639	7.0530579221286\\
5.41450735158116	7.08181397719676\\
5.50718774657607	7.12087271547296\\
5.59881398503458	7.161669298956\\
5.6902108755395	7.20385015057772\\
5.7811938096551	7.24696024967835\\
5.8715541083776	7.29058846296647\\
5.96175648312303	7.33470823571233\\
6.05208949922158	7.37926167542388\\
6.14207773297766	7.42384253767349\\
6.23185370510484	7.46835199637607\\
6.32199470246929	7.51289240400656\\
6.41188837859304	7.55693453074411\\
6.50203583406995	7.60048647405043\\
6.5924619931061	7.6434089065851\\
6.6832493074053	7.68564711186624\\
6.77440168812368	7.72701804874393\\
6.8659523476915	7.76761056849038\\
6.90489264391439	7.85978803620596\\
6.91647753064677	7.95917125137119\\
};
\addplot [color=mycolor3,line width=1.8pt,mark size=12.5pt,only marks,mark=o,mark options={solid},forget plot]
  table[row sep=crcr]{%
1.2007176764374	5.04250627212539\\
1.20092204394336	5.14360130196312\\
1.2023777160934	5.24475155962437\\
1.20628482128034	5.34536425874606\\
1.21389201796482	5.44643953999208\\
1.22616009115401	5.54571212815427\\
1.2442384021695	5.64407846567384\\
1.26943473253982	5.74209983216907\\
1.30231866779821	5.83769596650539\\
1.34297724343823	5.92924581757342\\
1.39224955191863	6.01752441469541\\
1.44908543697882	6.10032690954023\\
1.5132652775833	6.17779773069954\\
1.58371199430582	6.24934894437482\\
1.65995894260238	6.31535882188353\\
1.74065171782917	6.37556399973043\\
1.82523345916956	6.43047492482794\\
1.91222573017793	6.48003719419712\\
2.00208399602429	6.52528636366584\\
2.0933787883098	6.5661639397679\\
2.18672842379318	6.60352154909401\\
2.28182419916213	6.63766229587112\\
2.37822698535469	6.66519872898775\\
2.47702824580503	6.68784429694152\\
2.5758185882837	6.70808327032847\\
2.67482702206532	6.72646322662998\\
2.77472598628898	6.74342935489307\\
2.87486896687477	6.7591022606976\\
2.97471381856866	6.77359136082224\\
3.07543502487994	6.78720899563951\\
3.17503261599432	6.79980738792378\\
3.27480924335516	6.81165988642072\\
3.37447725341964	6.82281065212561\\
3.47539829708904	6.83346793022957\\
3.57568737591418	6.84348421282253\\
3.6767144468401	6.85304515533923\\
3.77664304354778	6.86202205695208\\
3.87683681067276	6.87058014386071\\
3.97705160186839	6.87872880770414\\
4.07704799121915	6.88647804272901\\
4.17806266526706	6.89394453086062\\
4.27833308073164	6.90101902345736\\
4.37902699966684	6.90780667813783\\
4.47923529327555	6.91433043059959\\
4.57954885319566	6.92147088942672\\
4.6799311226536	6.92987802441013\\
4.78025307209861	6.93858765250962\\
4.87611550716991	6.90885963736503\\
4.96717632336755	6.86552185901838\\
5.05780582679792	6.82211669923681\\
5.14820562836902	6.77877609785723\\
5.23856265552707	6.73546699307975\\
5.32908741233559	6.6921396626123\\
5.41998183202997	6.64875601849619\\
5.51075877233942	6.60562576051467\\
5.60169927918425	6.56272923873359\\
5.69251692655372	6.5203807916041\\
5.78374260573378	6.4786198702885\\
5.87556009566327	6.43778167584179\\
5.9678703484499	6.39834456381005\\
6.06052793279695	6.36056511584068\\
6.15401018991913	6.32417050553263\\
6.24789503868147	6.28918263580478\\
6.34213505498724	6.25552643446679\\
6.4369398312668	6.22326437187645\\
6.53184382823671	6.19061281717883\\
6.62503758540823	6.15432001676216\\
6.71573821115119	6.11213939357727\\
6.80186271189295	6.0610712069689\\
6.8769652021846	5.99504011128476\\
};

\addplot[area legend,solid,draw=black,fill=black,forget plot]
table[row sep=crcr] {%
x	y\\
2.32743746213213	9.04930358224796\\
1.8	9.2\\
1.8	9.2\\
1.8	9.2\\
1.8	9.2\\
1.8	9.2\\
2.85487492426427	8.89860716449592\\
2.8177067842983	8.82332241082026\\
3.2	8.8\\
2.89204306423024	8.97389191817159\\
2.85487492426427	8.89860716449592\\
}--cycle;
\node[right, align=left, text=black]
at (axis cs:2.7,9.2) {\Huge $\alpha_1$};

\addplot[area legend,solid,draw=black,fill=black,forget plot]
table[row sep=crcr] {%
x	y\\
1.56566820276498	7.4\\
1.1	7.4\\
1.1	7.4\\
1.1	7.4\\
1.1	7.4\\
1.1	7.4\\
2.03133640552995	7.4\\
2.03133640552995	7.31958061047956\\
2.4	7.4\\
2.03133640552995	7.48041938952044\\
2.03133640552995	7.4\\
}--cycle;
\node[right, align=left, text=black]
at (axis cs:1.5,7.65) {\Huge $\alpha_2$};

\addplot[area legend,solid,draw=black,fill=black,forget plot]
table[row sep=crcr] {%
x	y\\
1.8	5.30977212971078\\
1.8	4.95\\
1.8	4.95\\
1.8	4.95\\
1.8	4.95\\
1.8	4.95\\
1.8	5.66954425942156\\
1.69428741538846	5.66954425942156\\
1.8	5.95\\
1.90571258461154	5.66954425942156\\
1.8	5.66954425942156\\
}--cycle;
\node[right, align=left, text=black]
at (axis cs:2.05,5.45) {\Huge $\alpha_3$};
\end{axis}
\end{tikzpicture}%
\caption{Trajectory for $\phi_3$}\label{fig:3}
\end{subfigure}
\begin{subfigure}{0.48\textwidth}
%
%
\definecolor{mycolor1}{rgb}{0.00000,0.44700,0.74100}%
\definecolor{mycolor2}{rgb}{0.85000,0.32500,0.09800}%
\definecolor{mycolor3}{rgb}{0.92900,0.69400,0.12500}%
\begin{tikzpicture}[scale=0.45]

\begin{axis}[%
width=6.028in,
height=4.754in,
at={(1.011in,0.642in)},
scale only axis,
xmin=5,
xmax=10,
xlabel={x-position},
xmajorgrids,
ymin=4,
ymax=10,
ylabel={y-position},
ymajorgrids,
axis background/.style={fill=white},
title style={font=\LARGE},xlabel style={font=\LARGE},ylabel style={font=\LARGE},legend style={font=\LARGE},ticklabel style={font=\LARGE},
]
\addplot [color=mycolor1,line width=1.8pt,mark size=12.5pt,only marks,mark=o,mark options={solid},forget plot]
  table[row sep=crcr]{%
7.99840123197712	7.00043354799197\\
8.07354789938229	7.00019008411731\\
8.14861878968829	7.00038402644354\\
8.22362590966839	7.00039688353958\\
8.29870392474208	7.00040853416475\\
8.37371677758669	7.00041924270301\\
8.44871751045385	7.00042910881605\\
8.52373511371675	7.00043820852727\\
8.59874002609628	7.00044659912781\\
};
\addplot [color=mycolor2,line width=1.8pt,mark size=12.5pt,only marks,mark=o,mark options={solid},forget plot]
  table[row sep=crcr]{%
6.92990602676161	8.05778567683788\\
6.90063516485107	8.12777220287526\\
6.87890191155714	8.19994912109537\\
6.85781857910638	8.27209280724604\\
6.83445339246857	8.34336209511276\\
6.81099050431786	8.41470454049926\\
6.78754626137385	8.48602057018333\\
6.76410239957078	8.55737266521485\\
6.74067994123649	8.62870341352225\\
6.71728017855282	8.70001635631203\\
6.69391082269322	8.77129679162351\\
6.67055961336421	8.84259280840402\\
6.6472356860375	8.91388982330828\\
};
\addplot [color=mycolor3,line width=1.8pt,mark size=12.5pt,only marks,mark=o,mark options={solid},forget plot]
  table[row sep=crcr]{%
6.90064415073048	5.93908057405187\\
6.88746873832896	5.86516065167407\\
6.85952306313499	5.79543520421102\\
6.83087640876287	5.72596783020501\\
6.80645801546396	5.65503826670195\\
6.78214202891724	5.58395790072391\\
6.75782435846064	5.51293213285627\\
6.73349824118403	5.44195423025349\\
6.70914501525669	5.37098191557137\\
6.68474812858438	5.29998103230273\\
6.66033217818256	5.2290396118142\\
6.63587389849727	5.15810979316721\\
6.61137097420932	5.08720930315457\\
};

\addplot[area legend,solid,draw=black,fill=black,forget plot]
table[row sep=crcr] {%
x	y\\
8.20783410138249	7.4\\
7.9	7.4\\
7.9	7.4\\
7.9	7.4\\
7.9	7.4\\
7.9	7.4\\
8.51566820276498	7.4\\
8.51566820276498	7.31958061047956\\
8.7	7.4\\
8.51566820276498	7.48041938952044\\
8.51566820276498	7.4\\
}--cycle;
\node[right, align=left, text=black]
at (axis cs:8.2,7.6) {\Huge $\alpha_1$};

\addplot[area legend,solid,draw=black,fill=black,forget plot]
table[row sep=crcr] {%
x	y\\
7.16515958497107	8.5237841615195\\
7.3	8.1\\
7.3	8.1\\
7.3	8.1\\
7.3	8.1\\
7.3	8.1\\
7.03031916994213	8.94756832303901\\
6.98274445288636	8.91252706206059\\
6.95	9.2\\
7.0778938869979	8.98260958401744\\
7.03031916994213	8.94756832303901\\
}--cycle;
\node[right, align=left, text=black]
at (axis cs:7.225,8.75) {\Huge $\alpha_2$};

\addplot[area legend,solid,draw=black,fill=black,forget plot]
table[row sep=crcr] {%
x	y\\
6.51668784023658	5.85006352070973\\
6.65	6.25\\
6.65	6.25\\
6.65	6.25\\
6.65	6.25\\
6.65	6.25\\
6.38337568047315	5.45012704141946\\
6.33623530820772	5.48650178210488\\
6.3	5.2\\
6.43051605273858	5.41375230073404\\
6.38337568047315	5.45012704141946\\
}--cycle;
\node[right, align=left, text=black]
at (axis cs:6.225,5.975) {\Huge $\alpha_3$};
\end{axis}
\end{tikzpicture}%
\caption{Trajectory for $\phi_4$}\label{fig:4}
\end{subfigure}
\caption{Trajectories of the three robots.}
\end{figure*}

\begin{figure*}[t]
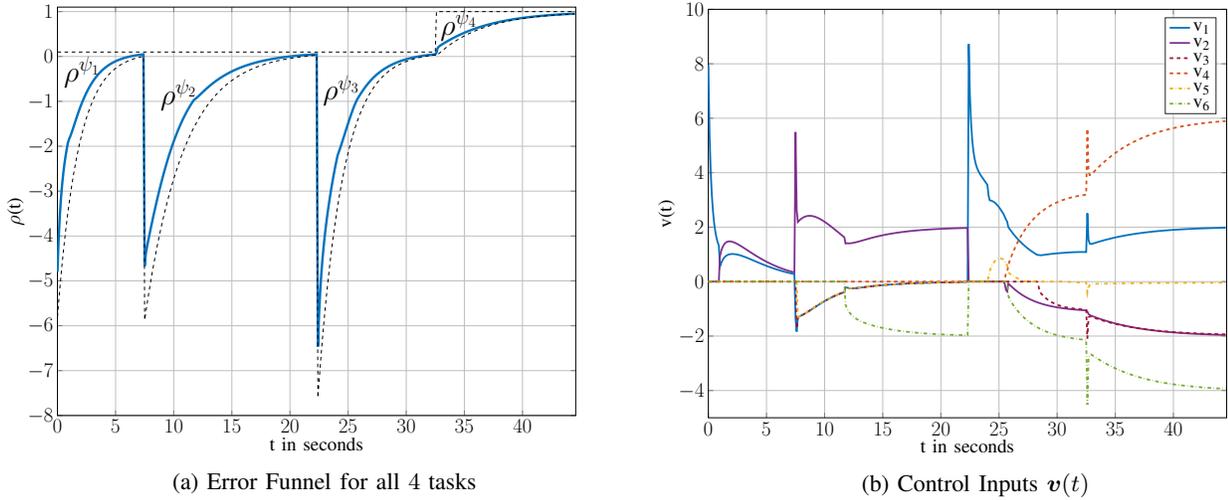

\centering
\begin{subfigure}{0.48\textwidth}
\input{figures/funnel}\caption{Error Funnel for all $4$ tasks}\label{fig:6}
\end{subfigure}
\begin{subfigure}{0.48\textwidth}
\input{figures/input}\caption{Control Inputs $\boldsymbol{v}(t)$}\label{fig:7}
\end{subfigure}
\caption{Time evolution of error and inputs}
\end{figure*}

We consider a multi-agent system with single integrator dynamics in $\mathbb{R}^2$ and deploy the well known consensus protocol \cite{mesbahi2010graph} with additional free inputs. The consensus protocol can be seen as the desire of the group to stay close to each other. In other words, assume $M$ agents where each agent $j\in\{1,\hdots,M\}$ is subject to the dynamics $\dot{\boldsymbol{x}}_j = \boldsymbol{v}_j$ with $\boldsymbol{x}_j\in \mathbb{R}^2$. The consensus protocol is then used as $\boldsymbol{v}_j = -\sum\limits_{k\in \mathcal{N}_j} (\boldsymbol{x}_j-\boldsymbol{x}_k) + \boldsymbol{u}_j$ where $\mathcal{N}_j$ denotes the neighborhood of the agent $j$. Using the graph Laplacian $L$ \cite{mesbahi2010graph}, we can express the dynamics as
\begin{align}\label{eq:dynamics1}
\dot{\boldsymbol{x}}(t)=-(L \otimes I_2) \boldsymbol{x}(t)+\boldsymbol{u}(t).
\end{align}
Comparing \eqref{eq:dynamics1} with \eqref{system_noise} reveals that $f(\boldsymbol{x})=-(L \otimes I_2) \boldsymbol{x}$ and $g(\boldsymbol{x})=I_M \otimes I_2=I_{2M}$, where $I_M$ is the $M\times M$ identity matrix. Note that Assumption \ref{assumption:1} is trivially satisfied. More specifically, assume three agents $\alpha_1$, $\alpha_2$, and $\alpha_3$ connected by means of a fixed and complete graph with a graph Laplacian $L=\begin{bmatrix}
1 & -1 & 0\\
-1 & 2 & -1\\
0 & -1 & 1
\end{bmatrix}$. Denote the robot position with $\boldsymbol{x}_j:=\begin{bmatrix} x_{j,1} & x_{j,2} \end{bmatrix}^T$ for $j\in\{1,2,3\}$. The initial positions are $\boldsymbol{x}_1(0):=\begin{bmatrix}1.1 & 3.1\end{bmatrix}^T$, $\boldsymbol{x}_2(0):=\begin{bmatrix}2 & 0.5\end{bmatrix}^T$, and $\boldsymbol{x}_3(0):=\begin{bmatrix}7 & 1.5\end{bmatrix}^T$. We also have five goal positions $A$, $B$, $C$, $D$, and $E$, which are located at $\boldsymbol{p}_A := \begin{bmatrix}6 & 4\end{bmatrix}^T$, $\boldsymbol{p}_B := \begin{bmatrix}1.2 & 9\end{bmatrix}^T$, $\boldsymbol{p}_C := \begin{bmatrix}1.2 & 7\end{bmatrix}^T$, $\boldsymbol{p}_D := \begin{bmatrix}1.2 & 5\end{bmatrix}^T$, and $\boldsymbol{p}_E := \begin{bmatrix}8 & 7\end{bmatrix}^T$. We use $\big(\|\boldsymbol{x}_j-\boldsymbol{p}_A\|_\infty<c\big) = \big(|x_{j,1}-p_{A,1}|<c\big) \wedge \big(|x_{j,2}-p_{A,2}|<c\big) = \big(x_{j,1}-p_{A,1}<c\big) \wedge \big(-x_{j,1}+p_{A,1}<c\big) \wedge \big(x_{j,2}-p_{A,2}<c\big) \wedge \big(-x_{j,2}+p_{A,2}<c\big)$ to ensure that $\|\boldsymbol{x}_j-\boldsymbol{p}_A\|_\infty=\max(|x_{j,1}-p_{A,1}|,|x_{j,2}-p_{A,2}|)<c$.

The robots are subject to the following sequential tasks: 1) Robot $\alpha_1$ moves to A within $7-10$ seconds. 2) Within the next $10-20$ seconds, $\alpha_1$, $\alpha_2$, and $\alpha_3$ move to $B$, $C$, and $D$, respectively. 3) $\alpha_1$ moves to $E$ within $5-15$ seconds. Additionally $\alpha_2$ and $\alpha_3$ form a triangular formation. 4) $\alpha_2$ and $\alpha_3$ always keep at least a distance of $1$ from $\alpha_1$ and disperse. More specifically, we have: $\theta:=F_{[7,10]}(\psi_1 \wedge F_{[10,20]}(\psi_2 \wedge F_{[5,15]}(\psi_3 \wedge \phi_4)))$ with $\psi_1:= \big(\|\boldsymbol{x}_1-\boldsymbol{p}_A\|_\infty<0.1\big)\wedge \psi_{Ass.3}$, $\psi_2:= \big(\|\boldsymbol{x}_1-\boldsymbol{p}_B\|_\infty<0.1\big) \wedge \big(\|\boldsymbol{x}_2-\boldsymbol{p}_C\|_\infty<0.1\big) \wedge \big(\|\boldsymbol{x}_3-\boldsymbol{p}_D\|_\infty<0.1\big)\wedge \psi_{Ass.3}$, $\psi_3:= \big(\|\boldsymbol{x}_1-\boldsymbol{p}_E\|_\infty<0.1\big) \wedge \big(1<x_{1,1}-x_{2,1}<1.2\big) \wedge \big(1<x_{1,1}-x_{3,1}<1.2\big) \wedge \big(1<x_{2,2}-x_{1,2}<1.2\big) \wedge \big(1<x_{1,2}-x_{3,2}<1.2\big)\wedge \psi_{Ass.3}$, and $\phi_4:=G_{[0,12]} \big((1<x_{1,1}-x_{2,1}) \wedge (1<x_{2,2}-x_{1,2}) \wedge (1<x_{1,1}-x_{3,1}) \wedge (1<x_{1,2}-x_{3,2})\wedge \psi_{Ass.3}\big)$ with $\psi_{Ass.3}:=(\|\boldsymbol{x}\|_\infty<100)$ to enforce Assumption \ref{assumption:4}.

The simulation result for all four tasks is displayed in Fig. \ref{fig:5}. In more detail, the trajectories for  $\phi_1$ and $\phi_2$ can be found in Fig. \ref{fig:1} and \ref{fig:2}, respectively. For $\phi_1$, the consensus dynamics bring the agents together, while at the same time the performance function $\gamma_1(t)$ forces $\alpha_1$ to approach and reach $A$, followed by agents $\alpha_2$ and $\alpha_3$. For the second task in Fig. \ref{fig:2}, each agent individually reaches its goals $B$, $C$, and $D$. The third task is shown in Fig. \ref{fig:3}, where we see that initially the robots gather and eventually form a triangular formation while $\alpha_1$ approaches $E$. In Fig. \ref{fig:4}, dispersion of the multi-agent system can be seen. To see that time bounds have been respected, Fig. \ref{fig:6} displays the different funnels. Fig. \ref{fig:7} shows that the control inputs are bounded and piecewise-continuous. To conclude, $\theta$ is satisfied with $\rs{\theta} > 0.05$. Note that due to the precision that we chose, e.g., $0.1$ in $\phi_1=F_{[7,10]} \big(\|\boldsymbol{x}_1-\boldsymbol{p}_A\|_\infty<0.1\big)$, $r$ can not exceed $0.1$. We remark that the control law is centralized and that simulations have been performed in real-time, which is possible due to the easy-to-implement feedback control law.

\section{Conclusion}
\label{sec:conclusion}
We considered nonlinear systems subject to a subset of signal temporal logic specifications. The imposed transient and steady-state behavior of the prescribed performance control approach was leveraged to satisfy atomic temporal formulas. A hybrid control strategy was then used to ensure that a finite set of atomic temporal formulas is satisfied.  A salient feature is that  the feedback control law is piecewise-continuous and robust with respect to disturbances and the specification, i.e.,  the specification is satisfied with a user-defined robustness.

Future work will include the extension of the derived methods to decentralized multi-agent systems with couplings in various forms. Local and global specifications will be subject of our work in this respect, as well as the feasibility of these coupled specifications. Furthermore, an extension of the expressivity, i.e., the signal temporal logic subset under consideration, will be investigated. %

\bibliographystyle{IEEEtran}
\bibliography{literature}

\addtolength{\textheight}{-12cm}   

\end{document}